\documentclass[12pt,a4paper]{article}
\usepackage[dvipdfmx,hiresbb]{graphicx}
\usepackage{latexsym}
\usepackage{amssymb}
\usepackage{amsmath}
\usepackage{amscd}
\usepackage{mathrsfs}
\usepackage{bm}

\newtheorem{Th}{Theorem}[section]
\newtheorem{Co}[Th]{Corollary}
\newtheorem{Lem}[Th]{Lemma}
\newtheorem{Exa}[Th]{Example}
\newtheorem{Rem}[Th]{Remark}

\newtheorem{Pro}[Th]{Proposition}
\newcommand{\demo}{\par\noindent{\it Proof. \/}\ }
\newcommand{\enD}{\hfill $\Box$ \vspace{3truemm}\par}

\newcommand{\R}{{\mathbb R}}
\newcommand{\lon}{\longrightarrow}

\begin{document}

\title{Geometric equivalence among smooth map germs } 

\author{Shyuichi IZUMIYA, Masatomo TAKAHASHI \\
and \\
Hiroshi TERAMOTO
\bigskip
\\
Dedicated to the memory of John N. MATHER}

\date{\today}

\maketitle

\begin{abstract} We consider equivalence relations among smooth map germs 
with respect to geometry of $G$-structures on the target space germ.
These equivalence relations are natural generalization of right-left equivalence (i.e., $\mathcal{A}$-equivalence)
in the sense of Thom-Mather depending on geometric structures on the target space germ.
Unfortunately,
these equivalence relations are not necessarily geometric subgroups in the sense of Damon (1984). 
However, we have interesting applications of these equivalence relations.
\end{abstract}
\renewcommand{\thefootnote}{\fnsymbol{footnote}}
\footnote[0]{2010 Mathematics Subject classification. Primary 58K40 ; Secondary 53C10}
\footnote[0]{Keywords. $G$-structure, $\mathcal{A}$-equivalence, Singularities of map germs}
\footnote[0]{This work is partially supported by Grant-in-Aid for Scientific Research (JSPS); (A) 26247006 (S.I),
(C) 17K05238 (M.T) and JST PRESTO; JPMJPR16E8 (H.T).} 

\section{Introduction}
{\ \ \ \,}
In the history of the theory of singularities of smooth mapping, the notion of
$\mathcal{A}$-equivalence (i.e. right-left equivalence or isomorphism) among smooth map germs in the sense of Mather is
the most natural equivalence (cf. \cite{MatherIII, MatherIV}) from the view point of differential topology.
In order to solve the stability problems of Thom \cite{TL}, Mather also introduced the notion of $\mathcal{K}$-equivalence, which played a key role in his theory.
Moreover, Tougeron \cite{Tougeron} introduced the notion of $\mathcal{K}[G]$-equivalence (it is $G$-equivalence 
in the terminology of Tougeron) for a linear Lie group $G$ which linearly acts on the target space germ.
If $G$ is a general linear group, then $\mathcal{K}[G]$-equivalence is $\mathcal{K}$-equivalence. Recently, there appeared several applications of $\mathcal{K}[G]$-equivalence (cf. \cite{ArnoldM,BruceM,Bruce-TariM,CPS,ITK,Ped,SF,TKITK})
which include applications to quantum physics etc.
In this paper we consider the case when the target space germ $(\R^{p},0)$ has a $G$-structure.
Then we introduce the notion of $\mathcal{A}[G]$-equivalence among smooth map germs analogous to $\mathcal{K}[G]$-equivalence.
If $G=GL(p,\R)$, then $\mathcal{A}[GL(p,\R)]$-equivalence is the original $\mathcal{A}$-equivalence in the sense of Mather.
If $G=\{I_{p}\}$ ($I_{p}$ is the unit matrix), then $\mathcal{A}[\{I_{p}\}]$-equivalence is $\mathcal{R}$-equivalence in
the sense of Mather \cite{MatherIII,MatherR}.
Therefore, $\mathcal{A}[G]$-equivalence is one of the direct and natural generalizations of
$\mathcal{A}$-equivalence.
Although $\mathcal{K}[G]$ is a geometric subgroup of $\mathcal{K}$ in the sense of Damon \cite{Damon}, $\mathcal{A}[G]$ is not necessarily a geometric subgroup of $\mathcal{A}$.
Thus the usual techniques of singularity theory cannot work generally.
Moreover, it is known that $\mathcal{A}$-equivalence implies $\mathcal{K}$-equivalence \cite{MatherIII}.
This fact does not hold for $\mathcal{A}[G]$ and $\mathcal{K}[G]$ generally.
The above properties are dependent on the Lie group $G.$
However, we can seek out the interesting examples of $\mathcal{A}[G]$-equivalence, which have been investigated recently (cf. \cite{CPS,D-R,Dufour,F-H,Fuku-Taka,Mancini,M-S,M-S2,West}).
Therefore it is worth while to study properties of $\mathcal{A}[G]$-equivalence for general Lie subgroup
$G\subset GL(p,\R).$ In this paper, we consider some fundamental properties of $\mathcal{A}[G]$-equivalence and give some interesting examples.
As a first step, we investigate the infinitesimal algebraic structure of $\mathcal{A}[G]$-equivalence.
\par
On the other hand, if we consider a $G'$-structure on the source space germ $(\R^{n},0),$
we also have the notion of $\mathcal{R}[G'],$ $\mathcal{A}[G';G]$ and $\mathcal{K}[G';G]$-equivalence
among smooth map germs, respectively. Even though there are interesting examples of those equivalence relations, we need longer pages for describing those equivalence relations, so that we only consider $\mathcal{A}[G]$-equivalence in this paper. 
\par
The organization of this paper is as follows.
In \S 2 we introduce $\mathcal{A}[G]$-equivalence and $\mathcal{R}\times G$-equivalence which are main subjects in this paper.
We briefly review algebraic properties of infinitesimal version of $\mathcal{A}$-equivalence following Mather \cite{MatherIII} in \S 3.
For the study of $\mathcal{A}[G]$-equivalence, we investigate algebraic properties of vector field
associated with a linear Lie group $G$ in \S 4.
Moreover, we calculate these vector fields for some examples of Lie groups.
Following the results of previous sections, we formulate the infinitesimal version of 
$\mathcal{A}[G]$-equivalence and investigate the relationship between $\mathcal{A}[G]$ and $\mathcal{R}\times G$
in \S 5.
In \S 6 we give several examples of $\mathcal{A}[G]$-equivalence, where
$G$ is $SO(p),$ $ SL(p,\R),$ $Sp(2p)$ or other cases.
There are other interesting cases which we do not mention here (for example, $G=SO_{0}(1,q)$, etc)
which will be investigated in elsewhere.
Following the observations in \S 6, we investigate relationships between $\mathcal{A}[G]$-equivalence and
$\mathcal{R}\times G$-equivalence in \S 7.
Finally we propose an important prospective problem.
\par
We assume that all map germs and manifolds are class $C^\infty$ unless stated otherwise.

\section{Geometric equivalence}
{\ \ \ \,}
We consider smooth map germs $f:(\R^n,0)\lon (\R^p,0).$
One of the most natural equivalence relations among map germs is {\it $\mathcal{A}$-equivalence} (i.e. an {\it isomorphism})
in the sense of Mather \cite{MatherIII}.
We say that smooth map germs $f,g:(\R^n,0)\lon (\R^p,0)$ are  {\it $\mathcal{A}$-equivalent} if there exist diffeomorphism germs
$\phi :(\R^n,0)\lon (\R^n,0)$ and $\psi :(\R^p,0)\lon (\R^p,0)$ such that
$\psi\circ f=g\circ \phi.$
We define the group of diffeomorphism germs on $(\R^{p},0):$ 
\[
{\rm Diff}\, (p)=\{\psi \ |\ \psi :(\R^{p},0)\lon (\R^{p},0):\ \mbox{diffeomorphism\ germ}\ \}.
\]
\par
In this paper we consider the case when the target space $\R^p$ has a geometric structure.
Let $G\subset GL(p,\R)$ be a linear Lie group. Then $G$ can be considered as a structure group of the tangent bundle of $\R^p,$ 
so that the group $G$ gives a $G$-structure on the target space $\R^p.$
We define natural geometric equivalence among smooth map germs with respect to $G$-structures,
which is a generalization of $\mathcal{A}$-equivalence as follows:
For a diffeomorphism germ $\psi :(\R^p,0)\lon (\R^p,0)$, we have the Jacobi matrix $J_\psi (y)$ 
at $y\in (\R^p,0).$
We say that smooth map germs $f,g:(\R^n,0)\lon (\R^p,0)$ are {\it $\mathcal{A}[G]$-equivalent} if 
there exist diffeomorphism germs $\phi:(\R^n,0)\lon (\R^n,0)$  and
 $\psi :(\R^p,0)\lon (\R^p,0)$ with $J_\psi (y)\in G$ for any $y\in (\R^p,0)$
 such that
 $
f\circ \phi =\psi \circ g.
$
This equivalence is not a geometric subgroup of $\mathcal{A}$ in the
sense of Damon \cite{Damon} generally.
The situation depends on the Lie group $G.$
 We consider the group of diffeomorphism germs with respect to $G$:
 \[
 {\rm Diff}[G](p)=\{\psi\in {\rm Diff}\, (p)\ | J_y\psi \in G\ \mbox{for\ any}\ y\in (\R^p,0)\}.
 \] 
 We remark that ${\rm Diff}[GL(p,\R)](p)={\rm Diff}\, (p)$ and ${\rm Diff}[\{I_{p}\}](p)=\{1_{\R^{p}}\}.$
 For $\psi \in {\rm Diff}[G](p),$ we say that $\psi$ is {\it isotopic to the identity}
 if there exists a family $\Psi:(\R^p\times \R,0\times [0,1])\lon (\R^{p},0)$ such that $\psi_t\in {\rm Diff}[G](p)$, $\psi_0=1_{\R^p}$
 and $\psi_1=\psi$,
 where $\psi_t(x)=\Psi (x,t).$
 We define
 \[
  {\rm Diff}_0[G](p)=\{\psi\in  {\rm Diff}[G](p)\ |\ \psi\ \mbox{is\ isotopic\ to\ the\ identity}\ \}.
 \]
We say that $f,g:(\R^n,0)\lon (\R^p,0)$ are $\mathcal{A}_0[G]$-equivalent if there exist a diffeomorphism germ  
$\phi :(\R^n,0)\lon (\R^n,0)$ and $\psi \in {\rm Diff}_0[G](p)$ such that $f\circ\phi =\psi\circ g.$
Following the definition of the left equivalence of Mather \cite{MatherIII}, we say that
$f,g :(\R^n,0)\lon (\R^p,0)$ are {\it $\mathcal{L}[G]$-equivalent} if there exists $\psi \in {\rm Diff}[G](p)$ such that
$f =\psi\circ g.$

Moreover, for any $A\in G\subset GL(p,\R),$ we have the natural linear isomorphism $\psi_A:(\R^p,0)\lon (\R^p,0)$ defined by 
$\psi_A(y)=A.^ty,$ where $^{t}y$ is the transposed column vector of $y=(y_{1},\dots ,y_{p})$ and $A.^ty$ is the matrix product.
We say that smooth map germs  $f,g:(\R^n,0)\lon (\R^p,0)$ are {\it $\mathcal{R}\times G$-equivalent} if
there exist a diffeomorphism germ $\phi:(\R^n,0)\lon (\R^n,0)$  and $A\in G$ such that
$
f\circ \phi =\psi_A \circ g.
$
If $G=\{I_{p}\}$, then both of  $\mathcal{A}[G]$-equivalence and $\mathcal{R}\times G$-equivalence are equal to $\mathcal{R}$-equivalence (i.e. right equivalence) in the sense of Mather \cite{MatherIII}. 
Moreover, $\mathcal{A}[GL(p,\R)]$-equivalence is $\mathcal{A}$-equivalence.
By definition, if $f,g$ are $\mathcal{R}\times G$-equivalent, then these are $\mathcal{A}[G]$-equivalent.
There are several interesting examples of these equivalence relations.
\par
We now define
$\mathcal{A}(n,p)={\rm Diff}\, (n)\times {\rm Diff}\,(p)$. Since $G\subset {\rm Diff}[G](p)\subset {\rm Diff}\,(p),$ we have subgroups 
$\mathcal{A}[G](n,p)={\rm Diff}\, (n)\times {\rm Diff}[G](p)$ and $(\mathcal{R}\times G)(n,p)={\rm Diff}\, (n)\times G$ of
$\mathcal{A}(n,p).$

\section{Infinitesimal structures of $\mathcal{A}$-equivalence}
{\ \ \ \, } Following \cite{MatherIII,MatherIV}, we briefly review the infinitesimal properties of $\mathcal{A}$-equivalence among map-germs. Let $\mathcal{E}_{n}$ be the local $\R$-algebra of function germs of $n$-variables at the
origin with the unique maximal ideal $\mathfrak{M}_{n}.$
For a map germ $f:(\R^n,0)\lon (\R^p,0),$ we have a pull-back $\R$-algebra homomorphism $f^*: \mathcal{E}_p\lon \mathcal{E}_n$
defined by $f^{*}(h)=h\circ f.$
We also consider an $\mathcal{E}_n$-module of germs of vector fields along $f,$ which is defined by
\[
\theta (f)=\left\{\sum_{j=1}^p \eta _j(x)\frac{\partial}{\partial y_j}\circ f\ \Bigm |\ \eta_j\in \mathcal{E}_n,\ (j=1,\dots ,p)\ \right\},
\]
where, $y=(y_1,\dots ,y_p)\in \R^p.$
Therefore, $\theta (f)$ is identified with 
\[
C^\infty (n,p)=\{h\ |\ h:(\R^n,0)\lon \R^p;\ \mbox{map\ germ}\}\cong \mathcal{E}_n\times \cdots \times \mathcal{E}_n=\mathcal{E}_n^p
\]
as 
$\mathcal{E}_{n}$-modules.
Moreover, $\mathfrak{M}_n\theta (f)=\mathfrak{M}_nC^\infty (n,p)=\mathfrak{M}_n^p$ is an $\mathcal{E}_n$-submodule of
$C^\infty (n,p)$ consisting of all map germs $(\R^n,0)\lon (\R^p,0)$. Therefore, we define the action of $\mathcal{A}(n,p)$ on
$\mathfrak{M}_nC^\infty (n,p)$ by
$\mu ((\phi,\psi),f)= \psi\circ f\circ \phi^{-1}.$ The orbit through $f$ is the set of all map germs which are $\mathcal{A}$-equivalent to $f$.
Since $\mathcal{A}[G](n,p)$ and $(\mathcal{R}\times G)(n,p)$ are subgroups of $\mathcal{A}(n,p),$ the above action induces the
actions of these subgroups on $\mathfrak{M}_nC^\infty (n,p).$ 
\par
We now consider formal tangent spaces of an $\mathcal{A}$-orbit.
The tangent space of $\mathfrak{M}_nC^\infty (n,p)$ at $f$ is defined to be the set of $d(c(t))/dt|_{t=0}$ for a curve $c(t)\in  \mathfrak{M}_nC^\infty (n,p)$
with $c(0)=f.$ We denote it by $T_{f}\mathfrak{M}_nC^\infty (n,p)$, which is called a (formal) {\it tangent space} of $\mathfrak{M}_nC^\infty (n,p)$ at $f.$
Since $c(t)(x)=f_t(x)$ with $f_0=f,$
we have
\[
T_{f}\mathfrak{M}_nC^\infty (n,p)=\left\{\sum_{j=1}^p \eta _j(x)\frac{\partial}{\partial y_j}\circ f\ \Bigm |\ \eta_j\in \mathfrak{M}_n (j=1,\dots ,p)\ \right\}.
\]
Therefore $T_{f}\mathfrak{M}_nC^\infty (n,p)$ can be identified with $\mathfrak{M}_n\theta (f).$
We also define an {\it extended tangent space} of $\mathfrak{M}_nC^\infty (n,p)$ at $f$ by  $T_{f}\mathfrak{M}_nC^\infty (n,p)_e=\theta (f).$
\par
Following Mather \cite{MatherIII}, 
we define a mapping $tf:\theta (n)\lon \theta(f)$ by $tf(\zeta)=df\circ\zeta,$
where $\theta(n)=\theta(1_{\R^n})$ and $df:T\R^n\lon T\R^p$ is the differential map of $f$.
We remark that $\theta (n)$ is the $\mathcal{E}_n$-module of vector field germs on $(\R^n,0).$
Then $tf$ is an $\mathcal{E}_n$-homomorphism.
We also define $\omega f:\theta (p)\lon \theta(f)$ by $\omega f(\xi)=\xi\circ f.$
Then $\theta (f)$ is an $\mathcal{E}_{p}$-module through the pull-back homomorphism $f^*:\mathcal{E}_{p}\lon \mathcal{E}_{n}.$
In this sense, $\omega f$ is an $\mathcal{E}_{p}$-homomorphism.
Therefore, $(\omega f,tf ,\theta (p),\theta (n),  \theta (f))$ is called a mixed homomorphism of finite type over 
$f^*:\mathcal{E}_{p}\lon \mathcal{E}_{n}$ in \cite{MatherIII}.
The notion of mixed homomorphisms plays a principal role in the theory of Mather in \cite{MatherIII,MatherIV}.
\par
In order to investigate $\mathcal{A}[G]$-equivalence, we need to investigate infinitesimal properties
of ${\rm Diff}[G](p).$

 \section{Algebraic structures of vector field germs with respect to $G$}
 {\ \ \ \,}
 In order to investigate general properties of the set of vector field germs with respect to a $G$-structure,
 we consider a linear Lie group $G\subset GL(p,\R)$ and the Lie algebra $\mathfrak{g}=T_IG\subset M_p(\R).$
 Here, $M_p(\R)$ is the Lie algebra of $p\times p$-matrices over $\R .$
\par
 For any $\psi\in {\rm Diff}[G](p),$ we define the {\it formal tangent space} of ${\rm Diff}[G](p)$ at $\psi$ by
  \[
 T_\psi {\rm Diff}[G](p)=\left\{ \frac{d\psi_t}{d t}|_{t=0}\Bigm |\ \psi_t\in {\rm Diff}[G](p)\ \mbox{for}\ t\in (\R,0), \psi_0=\psi\right\}.
 \]
 If $\psi =1_{\R^p},$ then 
 \[
 J_y\left(\frac{d\psi_t}{d t}|_{t=0}\right)=\frac{d(J_y\psi_t)}{dt}|_{t=0}\in T_IG=\mathfrak{g},
 \]
 for any $\psi_t\in {\rm Diff}[G](p)$ with $\psi_0=1_{\R^p}$ and any $y\in (\R^{p},0).$
  Since $d\psi_t/dt|_{t=0}$ is a vector field germ, we have
 \[
 \frac{d\psi_t}{d t}|_{t=0}(y)=\sum_{i=1}^p \eta _i(y)\frac{\partial}{\partial y_i},
 \]
 which can be identified with the map germ $\eta =(\eta_1,\dots ,\eta_p):(\R^p,0)\lon (\R^p,0).$
 Therefore, we have
 \[
 \left(\frac{\partial \eta_{i}}{\partial y_j}(y)\right)\in \mathfrak{g}\ \mbox{for\ any}\ y\in (\R^p,0).
 \]
 Since $\mathfrak{g}$ is a real vector space, $T_{1_{{\R^{p}}}} {\rm Diff}[G](p)$ is also a real vector space.
 Actually, we have
 \[
 T_{1_{{\R^{p}}}} {\rm Diff}[G](p)=\left\{\sum_{i=1}^p \eta _i(y)\frac{\partial}{\partial y_i}\Bigm | 
 \left(\frac{\partial \eta_{i}}{\partial y_j}(y)\right)\in \mathfrak{g}\ \mbox{for\ any}\ y\in (\R^p,0),\ \eta_i(0)=0\right\}.
\]
\par
 We now define $\R$-linear subspaces of $\theta (p)$ by
\[
\theta [G](p) =\left\{\sum_{i=1}^{p} \eta _i(y)\frac{\partial }{\partial y_{i}} \Bigm | \left(\frac{\partial \eta_{i}}{\partial y_j}(y)\right)\in \mathfrak{g}\ \mbox{for\ any}\ y\in (\R^p,0)\ \right\}
\]
and
\[
\theta [G]_0(p)=\left\{\sum_{i=1}^{p} \eta _i(y)\frac{\partial }{\partial y_{i}}\in \theta[G] (p) \Bigm | \eta_i(0)=0 \ \right\}.
\]
By definition, $\theta[G](p)$ and $\theta[G]_0(p)$ are $\R$-linear subspaces of
$\theta ( p).$ By the above arguments, we have $T_{1_{\R^p}}{\rm Diff}[G](p)=\theta [G]_0(p).$
For any $\psi\in {\rm Diff}_0[G](p)$, we consider $\psi_t\in {\rm Diff}[G](p)\ \mbox{for}\ t\in (\R,0)$
such that $\psi_0=1_{\R^p}.$
Then we define $\widetilde{\psi}_t=\psi\circ \psi_t\in {\rm Diff}[G](p),$ so that $\widetilde{\psi}_0=\psi.$
It follows that
\[
 d\psi \circ \frac{d\psi_t}{d t}|_{t=0}=\frac{d\widetilde{\psi}_t}{d t}|_{t=0}\in T_\psi {\rm Diff}[G](p),
 \]
 so that we have a linear isomorphism
 $t\psi :\theta [G]_0(p)\lon T_\psi {\rm Diff}[G](p)$ 
 defined by
 $t\psi(\eta)=d\psi\circ \eta.$
 Therefore, we can identify $\theta [G]_0(p)$ with $T_\psi {\rm Diff}[G](p)$ through $d\psi.$
Moreover, $\theta ( p)$ is an $\mathcal{E}_p$-module.
We also define an $\mathcal{E}_p$-module $\mathfrak{g}(\mathcal{E}_{p})$ by
\[
\mathfrak{g}(\mathcal{E}_{p})=\{\zeta\ |\ \zeta :(\R^{p},0)\lon \mathfrak{g}: C^\infty\ \}.
\]
\par
We now define a sub $\R$-algebra $\mathcal{E}_p[G]$ of $\mathcal{E}_p$ such that $\theta[G]_0(p)$ is a sub $\mathcal{E}_p[G]$-module of $\theta( p)$ as follows:
For $\lambda \in \mathcal{E}_p,$ we define a map germ
\[
{\rm grad}_y \lambda :\theta( p)\lon M_p(\mathcal{E}_p)
\]
by
\[
{\rm grad}_y \lambda (\eta)=\eta \otimes\left(\frac{\partial \lambda}{\partial y_1},\dots ,
\frac{\partial \lambda}{\partial y_p}
\right)=\begin{pmatrix}
\eta_1 \\
\vdots \\
\eta _p
\end{pmatrix}
\left(\frac{\partial \lambda}{\partial y_1},\dots ,
\frac{\partial \lambda}{\partial y_p}
\right)
=
\begin{pmatrix}
\eta_1\frac{\partial \lambda}{\partial y_1} &\cdots & \eta_1\frac{\partial \lambda}{\partial y_p}\\
\vdots & \vdots & \vdots \\
\eta_p\frac{\partial \lambda}{\partial y_1} & \cdots &\eta_p\frac{\partial \lambda}{\partial y_p}
\end{pmatrix},\]
where $\eta =(\eta_1,\dots ,\eta _p)$ and $y=(y_1,\dots ,y_p).$
Then we define a subset of $\mathcal{E}_p$ by
\[
\mathcal{E}_p[G]=\left\{ \lambda \in \mathcal{E}_p\ \Bigm |\ {\rm grad}_y \lambda(\eta)\in \mathfrak{g}(\mathcal{E}_p)\ \mbox{for}\
\eta\in \theta[G]_0(p)\right\}.
\]
Since $\mathfrak{g}(\mathcal{E}_p)$ is an $\mathcal{E}_p$-module,
\[
{\rm grad}_y(\lambda_1+\lambda_2)(\eta)={\rm grad}_y\lambda_1(\eta)+{\rm grad}_y\lambda_2(\eta)\in
\mathfrak{g}(\mathcal{E}_p)
\]
and
\[
{\rm grad}_y (\lambda_1\lambda_2)(\eta)=
\lambda _2{\rm grad}_y\lambda_1(\eta)+\lambda_1{\rm grad}_y \lambda_2(\eta)\in
\mathfrak{g}(\mathcal{E}_p)
\]
for any $\lambda _1,\lambda_2 \in \mathcal{E}_p[G]$
and  $\eta\in \theta[G]_0(p).$
This means that $\mathcal{E}_p[G]$ is an $\R$-algebra.
For any $\eta\in \theta[G]_0(p)$ and $\lambda \in \mathcal{E}_p$,
we have
\[
\left(\frac{\partial \lambda\eta_i}{\partial y_j}\right)={\rm grad}_y \lambda (\eta)+\lambda \left(\frac{\partial \eta_i}{\partial y_j}\right).
\]
It follows that $\lambda \eta\in \theta[G]_0(p)$ if and only if
${\rm grad}_y\lambda (\eta)\in \mathfrak{g}(\mathcal{E}_p),$
so that $\theta[G]_0$ is an $\mathcal{E}_p[G]$-module.
Then we have the following theorem:
\begin{Th}
Let $R$ be a sub $\R$-algebra of $\mathcal{E}_p$ such that $\theta[G]_0(p)$ is a sub $R$-module of $\theta( p).$ Then $R\subset \mathcal{E}_p[G].$
\end{Th}
\demo
We consider any $\lambda \in R\subset \mathcal{E}_p.$ For any $\eta\in \theta[G]_0(p),$
we have
\[
\left(\frac{\partial \lambda\eta_i}{\partial y_j}\right)={\rm grad}_y\lambda (\eta)+\lambda \left(\frac{\partial \eta_i}{\partial y_j}\right).
\]
Since $\eta, \lambda\eta\in \theta[G]_0(p),$ we have
$(\partial \eta_i/\partial y_j), (\partial \lambda\eta_i/\partial y_j)\in \mathfrak{g}(\mathcal{E}_p).$
It follows that ${\rm grad}_y\lambda (\eta)\in \mathfrak{g}(\mathcal{E}_p),$ so that
$\lambda\in \mathcal{E}_p[G].$
\enD
We call $\mathcal{E}_p[G]$ a {\it maximum $\R$-subalgebra of $\mathcal{E}_p$ with respect to
$\theta[G]_0(p).$}
Moreover, $\mathcal{E}_p[G]$ is a $C^\infty$-ring in the sense of \cite{Dubuc,Ishikawa}.
We say that an $\R$-subalgebra $R$ of $\mathcal{E}_p$ is a {\it $C^\infty$-subring} of $\mathcal{E}_p$
if 
\[
f(\lambda_1,\dots ,\lambda _r)\in R
\]
for any $\lambda _1,\dots \lambda _r\in R$ and $f\in \mathcal{E}_r.$
In this case $R$ is a local ring with the unique maximal ideal $\mathfrak{M}_R=\mathfrak{M}_p\cap R.$
\begin{Pro}
The maximum $\R$-subalgebra $\mathcal{E}_p[G]$ of $\mathcal{E}_p$ with respect to $\theta[G]_0(p)$
is a $C^\infty$-subring of $\mathcal{E}_p.$
\end{Pro}
\demo
For any $\lambda_1,\dots ,\lambda_r\in \mathcal{E}_p[G]$ and $f\in \mathcal{E}_r,$
we would like to show that 
$$f(\lambda_1,\dots ,\lambda _r)\in \mathcal{E}_p[G].$$ 
Since
\[
\frac{\partial f(\lambda_1(y),\dots ,\lambda _r(y))}{\partial y_i}=\sum_{j=1}^r \frac{\partial f(\lambda_1(y),\dots ,\lambda (y))}{\partial \lambda_j}\frac{\partial\lambda _j}{\partial y_i},
\]
we have 
\begin{eqnarray*}
{\rm grad}_y f(\lambda_1,\dots ,\lambda _r)(\eta)
&=&
\begin{pmatrix}
\eta_1\frac{\partial f}{\partial \lambda_1} &\cdots & \eta_1\frac{\partial f}{\partial \lambda_r}\\
\vdots & \ddots & \vdots  \\
\eta_p\frac{\partial f}{\partial \lambda_1} & \cdots & \eta_p\frac{\partial f}{\partial \lambda_r}
\end{pmatrix}
\begin{pmatrix}
\frac{\partial \lambda_1}{\partial y_1} & \cdots & \frac{\partial \lambda_1}{\partial y_p} \\
\vdots & \ddots &\vdots \\
\frac{\partial \lambda_r}{\partial y_1} & \cdots & \frac{\partial \lambda_r}{\partial y_p}
\end{pmatrix} \\
&=&\sum_{j=1}^r \frac{\partial f}{\partial \lambda _j}(\lambda(y)){\rm grad}_y(\lambda _j)(\eta),
\end{eqnarray*} 
for any $\eta =(\eta_1,\dots ,\eta_p)\in \theta[G]_0(p).$
By definition, ${\rm grad}_y(\lambda _j)(\eta)\in \mathfrak{g}(\mathcal{E}_p).$
Since $\mathfrak{g}(\mathcal{E}_p)$ is an $\mathcal{E}_p$-module, ${\rm grad}_y f(\lambda_1,\dots ,\lambda _r)(\eta)\in 
\mathfrak{g}(\mathcal{E}_p),$ for any $\eta =(\eta_1,\dots ,\eta_p)\in \theta([G])_0(p).$
This completes the proof.
\enD
\par
We say that $A$ is a {\it differentiable algebra} (or, {\it $DA$-algebra}) if $A$ is an $\R$-algebra and there exists a surjective algebra homomorphism
$\phi :\mathcal{E}_n\lon A$ for some $n\in \mathbb{N}.$
These algebras are local rings with maximal ideals denoted by $\mathfrak{M}_A.$
A homomorphism $\alpha :A\lon B$ of $DA$-algebras is an algebra homomorphism such that there exists a map germ $g:(\R^p,0)\lon (\R^n,0)$
and $\psi \circ g^*=\alpha\circ \phi$, where
$g^*:\mathcal{E}_n\lon \mathcal{E}_p$ is the pull-back homomorphism, $\phi:\mathcal{E}_n\lon A$ and $\psi:\mathcal{E}_p\lon B$ are
surjective homomorphisms as $\R$-algebras.
We say that $A\subset \mathcal{E}_n$ is a {\it $DA$-subalgebra} if $A$ is a $DA$-algebra and the inclusion $i:A\subset \mathcal{E}_n$ is a homomorphism of $DA$-algebras. This means that there is a map germ $\phi : (\R^n,0)\lon (\R^p,0)$ such that $\phi ^*(\mathcal{E}_p)=A.$
For DA-algebras, modules over $DA$-algebras and homomorphisms of $DA$-algebras, the Malgrange preparation theorem holds \cite{Malgrange}.
There exists a criterion when a $C^\infty$-subring is a $DA$-algebra \cite[Appendix]{Ishikawa2017}.
\begin{Pro}
Let $R\subset \mathcal{E}_p$ be a $C^\infty$-subring. Then $R$ is a $DA$-algebra if and only if
$R$ is finitely generated as $C^\infty$-ring.
\end{Pro}
\par
We now give some important examples.
\begin{Exa}\rm
We consider the special orthogonal group $SO( p)\subset GL(p,\R).$
A diffeomorphism germ $\psi :(\R^p,0)\lon (\R^p,0)$ with $J_\psi (x)\in SO(p)$ for any $x\in (\R^p,0)$
is an isometry germ.
In this case, the corresponding Lie algebra is
\[
\mathfrak{so}(p)=\{X\in M_{p}(\R)\ |\ ^{t}X=-X\}.
\]
For convenience, we consider the case when $p=2.$ In this case we have
\[
\theta [SO(2)](2)=\left\{\sum_{i=1}^{2} \eta _i(y)\frac{\partial }{\partial y_{i}}\ \Bigm |
\begin{pmatrix}
\frac{\partial \eta_{1}}{\partial y_{1}} & \frac{\partial \eta_{1}}{\partial y_{2}} \\
\frac{\partial \eta_{2}}{\partial y_{1}} & \frac{\partial \eta_{2}}{\partial y_{2}}
\end{pmatrix} 
\in \frak{so}(2)(\mathcal{E}_2) \right\}.
\]
Since $\frak{so}(2)$ is the Lie algebra of anti-symmetric matrices,
we have
\[
\frac{\partial \eta_{1}}{\partial y_{1}}(y_1,y_2)=\frac{\partial \eta_{2}}{\partial y_{2}}(y_1,y_2)=0,\ \frac{\partial \eta_{1}}{\partial y_{2}}(y_1,y_2)=-\frac{\partial \eta_{2}}{\partial y_{1}}(y_1,y_2)
.\]
It follows that $\eta _1(y_1,y_2)=\eta _1(y_2)$ and $\eta _2(y_1,y_2)=\eta _2(y_1).$
Therefore, we have $\eta _1(y_2)=a_1+b_1y_2+\xi _1(y_2)$ and $\eta _2(y_1)=a_2+b_2y_1+\xi _2(y_1)$
for some $a_i,b_i\in \R$ and $\xi_i\in \mathfrak{M}^2_1.$
Thus, $(\partial \eta _1/\partial y_2)(y_1,y_2)=-(\partial \eta _2/\partial y_1)(y_1,y_2)$ means that $b_1=-b_2$ and $(d\xi_1/dy_2)(y_2)
=-(d\xi_2/dy_1)(y_1).$
The last equality means that $\xi_1(y_2)=\xi_2(y_1)=0.$
Hence, 
\[
\begin{pmatrix}
\frac{\partial \eta_{1}}{\partial y_{1}} & \frac{\partial \eta_{1}}{\partial y_{2}} \\
\frac{\partial \eta_{2}}{\partial y_{1}} & \frac{\partial \eta_{2}}{\partial y_{2}}
\end{pmatrix}\in \frak{so}(2)(\mathcal{E}_2)\ \mbox{if\ and\ only\ if}\
\eta_{1}(y_{2})=a_{1}+by_{2},\eta_{2}(y_{1})=a_{2}-by_{1}\ \mbox{for}\ a_{i},b\in \R.
\]
It follows that 
\[
\theta (2)\supset \theta [SO(2)](2)=\left\{a_{1}\frac{\partial }{\partial y_{1}}+a_{2}\frac{\partial }{\partial y_{2}}+b\left(y_{2}\frac{\partial }{\partial y_{1}}-y_{1}\frac{\partial }{\partial y_{2}}\right)\
\Bigm | a_{i},b\in \R\right\}.
\]
By the similar arguments to above, we have
\[
\theta [SO(p)](p)=\left\langle \frac{\partial }{\partial y_{i}}\Bigm | i=1,\dots ,p\right\rangle_{\R}
+\left\langle y_{j}\frac{\partial }{\partial y_{i}}-y_{i}\frac{\partial }{\partial y_{j}}\Bigm | 1\leq i< j\leq p\right\rangle_{\R}.
\]
It follows that
\[
\theta [SO(p)]_0(p)=\left\langle y_{j}\frac{\partial }{\partial y_{i}}-y_{i}\frac{\partial }{\partial y_{j}}\Bigm | 1\leq i< j\leq p\right\rangle_{\R}.
\]
\par
We can easily show that 
$$
\left\{y_{j}\frac{\partial }{\partial y_{i}}-y_{i}\frac{\partial }{\partial y_{j}}\Bigm | 1\leq i< j\leq p
\right\}
$$
is linearly independent, so that $\dim_{\R}\theta [SO(p)]_0(p)=p(p-1)/2=\dim_{\R}\mathfrak{so}(p).$
For any $\lambda (y)=\sum_{1\leq i< j\leq p} \lambda_{ij}(y_{j}\frac{\partial }{\partial y_{i}}-y_{i}\frac{\partial }{\partial y_{j}}),$ we have
\[
\Lambda =\begin{pmatrix}
0 & \lambda_{12}& \lambda_{13} &\cdots & \lambda _{1p} \\
-\lambda_{12} & 0 & \lambda_{23}&\cdots & \lambda_{2p}\\
\vdots & \vdots &\vdots & \ddots & \vdots \\
-\lambda_{1p} & -\lambda_{2p} & -\lambda_{3p}&\cdots & 0
\end{pmatrix}
\in \mathfrak{so}(p).
\]
It follows that we have the following proposition.
\begin{Pro} With the above notations, we have
$\theta[SO(p)]_{0}(p)\cong \mathfrak{so}(p)$ as Lie algebras.
\end{Pro}
Moreover, we have the following theorem.
\begin {Th} With the above notations, we have
${\rm Diff}_0\, [SO(p)] (p)=SO(p).$
\end{Th}
\demo
Since $SO(p)$ is connected,
we always have ${\rm Diff}_0\, [SO(p)](p)\supset SO(p).$
For any $\phi \in {\rm Diff}_0 \, [SO(p)](p),$ we consider $\phi_t\in {\rm Diff}\, [SO(p)](p)$ such that $\phi_0=1_{\R^p}$
and $\phi_1=\phi.$
Then, for any $t_0\in (\R,[0, 1]),$ we have
\[
 \frac{d\phi_t}{d t}|_{t=t_0}(x)=\sum_{i=1}^p \eta _i(y)\frac{\partial}{\partial y_i}\circ \phi_{t_0}\ \mbox{such\ that}\ \left(\frac{\partial \eta_{i}}{\partial y_j}(y)\right)\in \mathfrak{so}(p)\ \mbox{and}\ \eta_{i}(0)=0,
\]
for any $y\in (\R^p,0).$ By the previous arguments, $\eta _i(y)$ are linear function germs.
Since $\phi_t(0)=0,$ we can write
$\phi_t(y)=A(t).y+h_t(y)$ for some map germ $A:(\R,[0,1])\lon GL(p,\R)$ and $h_t\in \mathfrak{M}_p^2C^\infty (p,p).$
Thus
\[
 \frac{d\phi_t}{d t}|_{t=t_0}(y)=\frac{d A}{dt}|_{t=t_0}x+\frac{dh_t}{dt}|_{t=t_0}(y),
 \]
 so that $(dh_t/dt)|_{t=t_0}(y)=0.$
 This means that $h_t(y)=h_0(y),$ so that $\phi_t(y)=A(t).y+h_{0}(y).$
 However, $1_{\R^p}(y)=\phi_0(y)=A(0).y+h_{0}(y),$ so that $h_{0}(y)=0.$
Therefore, $\phi(y)=A(1).y.$ It follows that $A(1)=J_{\phi}(y)\in SO(p).$ 
 Thus, we have $\phi\in SO(p).$
\enD
\par
We denote ${\rm Diff}^\omega\, [SO(p)] (p)$ as the subgroup of ${\rm Diff}\, [SO(p)](p)$ consisting of real analytic germs of diffeomorphisms.
By the above arguments in the proof, if we have $\phi (y)=A.y+h(y)\in {\rm Diff}\, [SO(p)] (p)$ for some  $A\in SO(p)$ and $h(y)\in \mathfrak{M}_p^\infty C^\infty (p,p)$ with $h\not= 0,$
then $\phi$ is not isotopic to the identity.
Therefore, $\phi\notin SO(p).$
This proves that ${\rm Diff}^\omega\, [SO(p)] (p)=SO(p).$
We remark that we do not know ${\rm Diff}\, [SO(p)](p)=SO(p)$ or not.
\par
We now consider $\mathcal{E}_p[SO(p)]$. By definition, $\lambda \in \mathcal{E}_p[SO(p)]$ if and only if
${\rm grad}_y\lambda (\eta)\in \frak{so}(p)(\mathcal{E}_p)$ for any $\eta\in \theta [SO(p)]_0(p).$
The last condition is equivalent to
$y_i(\partial \lambda/\partial y_j)=-y_j(\partial \lambda/\partial y_i)$ for $i,j=1,\dots ,p.$
In particular, $y_i(\partial \lambda/\partial y_i)=0$ as function germs for $i=1,\dots ,p.$
It follows that $\lambda $ is a constant germ.
This means that $\mathcal{E}_p[SO(p)]=\R$. Since $\theta [SO(p)](p)$ is a finite dimensional $\R$-vector space, it is 
a finitely generated $\mathcal{E}_p[SO(p)]$-module.
\end{Exa}
\begin{Exa}\rm
We consider the special linear group $SL( p,\R)\subset GL(p,\R).$
In this case $\mathfrak{sl}(p,\R)$ is the Lie algebra of traceless $p \times p$-real matrices:
\[
\mathfrak{sl}(p,\R)=\{X\in M_p(\R)\ |\ {\rm Trace}\, X=0\}.
\]
A diffeomorphism germ $\psi :(\R^p,0)\lon (\R^p,0)$ with $J_\psi (y)\in SL(p,\R)$ for any $y\in (\R^p,0)$
is a volume preserving diffeomorphism germ.
Therefore, the algebraic structure of $\theta [SL(p,\R)](p)$ might be deeply related to
that of the space of differential forms.
Actually, we have the following proposition.
\begin{Pro}
Let $d(\Omega^{(p-2)})$ be the vector space of germs of exact differential $(p-1)$-forms. Then
$d(\Omega^{(p-2)})$ and $\theta[SL(p,\R)](p)$ are isomorphic as $\R$-vector spaces. Hence,
$\theta[SL(p,\R)](p)$ is an infinite dimensional vector space.
\end{Pro}
\demo
By definition, we have
\[
\theta [SL(p,\R)](p)=\left\{\eta =\sum_{i=1}^p \eta_i\frac{\partial }{\partial y_i}
\Bigm | \left(\frac{\partial \eta_i}{\partial y_j}\right)(y)\in \mathfrak{sl}(p,\R)\right\}.
\]
For any $\eta=\sum_{i=1}^p \eta _i(y)\frac{\partial }{\partial y_i}\in \theta (p),$
we define
\[
\omega _\eta =\sum_{k=1}^p (-1)^{k-1}\eta _k(y)dy_1\wedge \cdots \wedge dy_{k-1}\wedge dy_{k+1}\wedge\cdots \wedge dy_p.
\]
Then we have
\[
d\omega_\eta =\sum_{k=1}^p \frac{\partial \eta _k}{\partial y_k}(y)dy_1\wedge\cdots \wedge dy_p={\rm div}\, (\eta)dy_1\wedge\cdots \wedge dy_p.
\]
Therefore, $\eta\in \theta [SL(p,\R)](p)$ if and only if $d\omega_\eta =0.$
Then we can define a mapping
$\Phi : \theta [SL(p,\R)](p)\lon d(\Omega^{(p-2)})$ by
$\Phi (\eta)=\omega_{\eta}.$
For any $\omega\in d(\Omega^{(p-2)})$, by the Poincar\'e lemma, there exists a germ of $(p-2)$-form $\theta\in \Omega ^{p-2}$ such that $d\theta =\omega.$
Since $d\theta =\omega$ is a germ of $(p-1)$-form, it is written by
\[
d\theta=\sum_{k=1}^p \zeta _k(y)dy_1\wedge \cdots \wedge dy_{k-1}\wedge dy_{k+1}\wedge\cdots \wedge dy_p.
\]
If we define $\eta_k(y)=(-1)^k \zeta _k(y)$ and $\eta=\sum_{i=1}^p \eta _i(y)\frac{\partial }{\partial y_i},$
then $\omega _\eta =d\theta=\omega.$
Since $dd\theta =0,$ we have ${\rm div}\, (\eta)=0,$ so that $\eta \in \theta [SL(p,\R)](p).$
This means that $\Phi (\eta)=\omega.$
By definition, $\omega_{\eta_1+\eta_2}=\omega_{\eta_1}+\omega _{\eta _2}$ and $\omega_{c\eta}=c\omega_{\eta}$ for any $c\in \R.$ 
We also have that $\omega_\eta=0$ if and only if $\eta =0.$
Thus $\Phi$ is an $\R$-linear isomorphism.
\enD
\par
In order to simplify the arguments, we consider the case $p=2.$
For a vector field 
\[
\eta =\eta _1(y_1,y_2)\frac{\partial}{\partial y_1}+\eta _2(y_1,y_2)\frac{\partial}{\partial y_2},
\]
$\eta \in \theta [SL(2,\R)](2)$ if and only if $\partial \eta_1/\partial y_1+\partial \eta _2/\partial y_2=0.$
Therefore, we have
\[
d(\eta_2dy_1-\eta_1dy_2)=\left(\frac{\partial \eta_1}{\partial y_1}+\frac{\partial \eta_2}{\partial y_2}\right)
dy_2\wedge dy_1=0.
\]
By the Poincar\'e lemma, there exists a function germ $f\in \mathcal{E}_2$ such that
\[
\eta_1=-\frac{\partial f}{\partial y_2},\ \eta _2=\frac{\partial f}{\partial y_1}.
\]
Here, we can choose $f\in \mathcal{E}_{2}$ with $f(0)=0.$
Thus we have
\[
 \theta [SL(2,\R)](2)=\left\{-\frac{\partial f}{\partial y_2}\frac{\partial }{\partial y_1}
 +\frac{\partial f}{\partial y_1}\frac{\partial }{\partial y_2}\ \Bigm | f\in \mathfrak{M}_2\ \right\}.
 \]
 We define a mapping $\Delta :\mathfrak{M}_2\lon  \theta [SL(2,\R)](2)$
 by 
 \[
 \Delta (f)=-\frac{\partial f}{\partial y_2}\frac{\partial }{\partial y_1}
 +\frac{\partial f}{\partial y_1}\frac{\partial }{\partial y_2}.
 \]
 Then $\Delta $ is an $\R$-linear isomorphism.
 It follows that
 $ \theta [SL(2,\R)](2)\cong \mathfrak{M}_2$ as $\R$-vector spaces.
 Therefore, $\theta [SL(2,\R)](2)$ is an infinite dimensional $\R$-vector space.
 If $\eta_{1}(0)=\eta _{2}(0)=0,$ then $f\in \mathfrak{M}^{2}_{2},$ so that
 $ \theta [SL(2,\R)]_{0}(2)\cong \mathfrak{M}^{2}_2.$
 \par
 On the other hand, suppose that $\lambda (y_1,y_2)\in \mathcal{E}_2[SL(2,\R)].$
Then we have 
$$
\eta_1(y_1,y_2)\frac{\partial \lambda}{\partial y_1}+\eta_2(y_1,y_2)\frac{\partial \lambda}{\partial y_2}=0,
$$
where $\partial \eta_1/\partial y_1+\partial \eta _2/\partial y_2=0$ and $\eta_1(0,0)=\eta_2(0,0)=0.$
If we choose linear function germs $\eta _1=a_1y_1+a_2y_2, \eta_2=b_1y_1+b_2y_2,$ then $b_2=-a_1.$
It follows that
\[
(a_1y_1+a_2y_2)\frac{\partial \lambda}{\partial y_1}+(b_1y_1-a_1y_2)\frac{\partial \lambda}{\partial y_2}=0.
\]
If we substitute $a_1=0,a_2=b_1=1,$ then
\[
y_2\frac{\partial \lambda}{\partial y_1}+y_1\frac{\partial \lambda}{\partial y_2}=0.
\]
Moreover, if we substitute $a_1=1,a_2=b_1=0,$ then
\[
y_1\frac{\partial \lambda}{\partial y_1}-y_2\frac{\partial \lambda}{\partial y_2}=0.
\]
Therefore, we have a system of linear equations:
\[
\left\{
\begin{matrix}
y_2\frac{\partial \lambda}{\partial y_1}+y_1\frac{\partial \lambda}{\partial y_2}=0 \\
y_1\frac{\partial \lambda}{\partial y_1}-y_2\frac{\partial \lambda}{\partial y_2}=0
\end{matrix}
\right.
\]
If $y_1^2+y_2^2\not= 0,$ we have $\partial\lambda/\partial y_1 =\partial \lambda/\partial y_2=0.$
Taking the limit $(y_1,y_2)\lon (0,0)$, we also have 
$(\partial\lambda/\partial y_1)(0,0) =(\partial \lambda/\partial y_2)(0,0)=0.$
Therefore, $\lambda$ is a constant function.
Hence, we have $\mathcal{E}_2[SL(2,\R)]=\R.$
Therefore, $ \theta [SL(p,\R)](2)$ and $\theta[SL(2,\R)]_{0}(2)$ are not finitely generated $\mathcal{E}_2[SL(2,\R)]$-modules. 
\end{Exa}
\begin{Exa}\rm
We consider the symplectic linear group $Sp(2p,\R)\subset GL(2p,\R),$ which is defined by
\[
Sp(2p,\R)=\left\{A\in GL(2p,\R)\ |\ ^t\!AJ_{2p}A=J_{2p}\right\}.
\]
Here,
\[
J_{2p}=\begin{pmatrix}
0 & I_p \\
-I_p & 0
\end{pmatrix}.
\]
In this case, the corresponding Lie algebra is
\[
\frak{sp}(2p,\R)=\{X\in M_{2p}(\R)\ |\ ^t\!XJ_{2p}+J_{2p}X=O\}.
\]
A diffeomorphism germ $\psi :(\R^{2p},0)\lon (\R^{2p},0)$ with $J_\psi (y)\in Sp(2p,\R)$ for any $y\in (\R^{2p},0)$
is a symplectic diffeomorphism germ for the canonical symplectic structure $\omega$ on $\R^{2p}.$
\par
We now consider the case $p=1.$ In this case it is easy to show that $Sp(2,\R)=SL(2,\R),$
so that $\mathfrak{sp}(2,\R)=\mathfrak{sl}(2,\R)$ and
\[
\theta [Sp(2,\R)](2)=\theta [SL(2,\R)](2)=\left\{-\frac{\partial f}{\partial y_2}\frac{\partial }{\partial y_1}
 +\frac{\partial f}{\partial y_1}\frac{\partial }{\partial y_2}\ \Bigm | f\in \mathfrak{M}_2\ \right\}.
\]
On the other hand, $\mathcal{E}_2[Sp(2,\R)]=\mathcal{E}_2[SL(2,\R)]=\R.$
For $p\geq 2,$ the structure of $\theta [Sp(2p,\R)]$ is complicated.
\end{Exa} 
 \begin{Exa}\rm 1) We consider the following example:
 $$
 D^{*}(p_{1},p_{2})=\left\{ \begin{pmatrix}
 A & O \\
 O & B
 \end{pmatrix}
 \in GL(p,\R) \Bigm | A\in GL(p_1,\R),\ B\in GL(p_2,\R)\right\},
 $$
where $p=p_1+p_2.$ Actually, we have $D^{*}(p_{1},p_{2})=GL(p_{1},\R)\oplus GL(p_{2},\R).$
 Then we have two Lie subgroups 
 \begin{eqnarray*}
H&=& \left\{ \begin{pmatrix}
 A & O \\
 O & I_{p_{2}}
 \end{pmatrix}
 \bigm | A\in GL(p_1,\R)\right\}=GL(p_{1},\R)\oplus\{I_{p_{2}}\},\\
 K&=&\left\{ \begin{pmatrix}
 I_{p_{1}} & O \\
 O & B
 \end{pmatrix}
 \bigm | B\in GL(p_2,\R)\right\}=\{I_{p_{1}}\}\oplus GL(p_{2},\R) .
 \end{eqnarray*}
 In this case, we have $D^{*}(p_{1},p_{2})=HK\cong H\times K$. The corresponding Lie algebras are
 \[
\mathfrak{h}=
 \left\{ \begin{pmatrix}
 A & O \\
 O & O
 \end{pmatrix}
\in M_p(\R) \Bigm | A\in M_{p_1}(\R)\right\}=M_{p_{1}}(\R)\oplus\{O\},
\]
\[
\mathfrak{k}=
 \left\{ \begin{pmatrix}
 O & O \\
 O & B
 \end{pmatrix}
\in M_p(\R) \Bigm | B\in M_{p_2}(\R)\right\}=\{O\}\oplus M_{p_{2}}(\R)
\]
and
$\mathfrak{d}(p_{1},p_{2})=\mathfrak{h}\oplus\mathfrak{k}=M_{p_{1}}(\R)\oplus M_{p_{2}}(\R).$
Moreover, we have 
 \[
\theta [H]_0(p)=\frak{M}_{p_1}\theta (p_1),\ \theta [K]_0(p)=\frak{M}_{p_2}\theta (p_2)
 \]
 and
 \[
 \theta[D^{*}(p_{1},p_{2})]_0(p)=\theta[H]_0(p)\oplus \theta[K]_0(p)\cong \frak{M}_{p_1}\theta (p_1)\oplus \frak{M}_{p_2}\theta (p_2)
 \]
 as $\R$-vector spaces.
 We can show that $\mathcal{E}_p[D^{*}(p_{1},p_{2})]=\R$, $\mathcal{E}_p[H]=\{\lambda (y)\ |\ y=(y_1,\dots y_{p_1})\in (\R^{p_1},0)\}=\mathcal{E}_{p_{1}}$
 and $\mathcal{E}_p[K]=\{\lambda (y)\ |\ y=(y_{p_1+1},\dots y_{p})\in (\R^{p_{2}},0)\}=\mathcal{E}_{p_{2}},$
 where $\R^p=\R^{p_1}\times\R^{p_2}.$
 It follows that 
$
\mathcal{E}_p[D^{*}(p_{1},p_{2})]=\mathcal{E}_p[H]\cap \mathcal{E}_p[K]=\R.
$
By definition, all of the above rings are DA-algebras.
 \par\noindent
 2) We consider the following example:
 $$
 T^{*}_{r}(p_{1},p_{2})=\left\{ \begin{pmatrix}
 A & B \\
 O & C
 \end{pmatrix}
 \in GL(p,\R) \Bigm | A\in GL(p_1,\R),\ C\in GL(p_2,\R)\right\}.
 $$
 Here, we write $T^{*}_{r}(p_{1},p_{2})=GL(p_{1},\R)\widetilde{\oplus}_{r} GL(p_{2},\R).$
Then we have two Lie subgroups 
 \begin{eqnarray*}
 N&=& \left\{ \begin{pmatrix}
 A & B \\
 O & I_{p_{2}}
 \end{pmatrix}
 \bigm | A\in GL(p_1,\R)\right\}=GL(p_{1},\R)\widetilde{\oplus}_{r}\{I_{p_{2}}\},\\ 
 K&=&\left\{ \begin{pmatrix}
 I_{p_{1}} & O \\
 O & C
 \end{pmatrix}
 \bigm | C\in GL(p_2,\R)\right\}=\{I_{p_{1}}\}\oplus GL(p_{2},\R) . 
\end{eqnarray*}
In this case, $N$ is a normal subgroup of $T^{*}_{r}(p_{1},p_{2})$. Then we have $T^{*}_{r}(p_{1},p_{2})=NK\cong N\rtimes K$ (i.e. the semi-direct product).
The corresponding Lie algebras are 
 \[
 \mathfrak{n}=
 \left\{ \begin{pmatrix}
 A & B \\
 O & O
 \end{pmatrix}
\in M_p(\R) \Bigm | (A,B)\in M_{p_1\times p}(\R)\right\},
\]
\[
 \mathfrak{k}=
 \left\{ \begin{pmatrix}
 O & O \\
 O & C
 \end{pmatrix}
\in M_p(\R) \Bigm | C\in M_{p_2}(\R)\right\}=\{O\}\oplus M_{p_{2}}(\R)
\]
and
$\mathfrak{t}_{r}(p_{1},p_{2})=\mathfrak{n}\oplus\mathfrak{k}.$
Moreover, we have 
 \[
 \theta[N]_0(p)=\frak{M}_p\theta (\pi_{p_1}),\ \theta[K]_0(p)=\frak{M}_{p_2}\theta (p_2)
 \]
 and 
 \[
 \theta[T^{*}_{r}(p_{1},p_{2})]_0(p)=\theta[N]_0(p)\oplus \theta[K]_0(p)=\frak{M}_p\theta (\pi_{p_1})\oplus \frak{M}_{p_2}\theta (p_2)
 \]
 as an $\R$-vector space, where $\pi_{p_1}:\R^p=\R^{p_1}\times\R^{p_2}\lon \R^{p_1}$ is the canonical projection
 and $\theta (\pi_{p_{1}})=C^{\infty}(p,p_{1}).$
 We can show that $\mathcal{E}_p[T^{*}_{r}(p_{1},p_{2})]= \mathcal{E}_p[K]=\{\lambda (y)\ |\ y=(y_{p_1+1},\dots y_{p})\in (\R^{p_2},0)\}=\mathcal{E}_{p_2}$, 
 $\mathcal{E}_p[N]=\mathcal{E}_p,$
 where $\R^p=\R^{p_1}\times\R^{p_2}.$
Since $\mathcal{E}_p[T^{*}_{r}(p_{1},p_{2})]=\mathcal{E}_p[K]\subset \mathcal{E}_p=\mathcal{E}_p[N]$, we have $\mathcal{E}_p[T^{*}_{r}(p_{1},p_{2})]=\mathcal{E}_p[N]\cap \mathcal{E}_p[K].$
Therefore, $\theta[T^{*}_{r}(p_{1},p_{2})]_0(p)$ is an $\mathcal{E}_{p_2}$-module, which is not finitely generated.
However, $\theta[N]_0(p)$ is a finitely generated  $\mathcal{E}_p[N]$-module and $\theta[K]_0(p)$ is a finitely generated  $\mathcal{E}_p[K]$-module,
respectively.
By definition, all of the above rings are DA-algebras.
\end{Exa}
\par
Moreover, we have the following examples.
\begin{Exa}\rm
1) We consider the subgroup $SO(p)\cap T^{*}_{r}(p_{1},p_{2})$
of $T^{*}_{r}(p_{1},p_{2})$ for $p=p_{1}+p_{2}.$ Then we can show that
\[
SO(p)\cap T^{*}_{r}(p_{1},p_{2})=\left\{ \begin{pmatrix}
 A & O \\
 O & B
 \end{pmatrix}
  \Bigm | A\in SO(p_1),\ B\in SO(p_2)\right\}=SO(p_{1})\oplus SO(p_{2}).
  \]
It follows that
\[
\theta [SO(p)\cap T^{*}_{r}(p_{1},p_{2})]_{0}(p)\cong \mathfrak{so}(p_{1})\oplus \mathfrak{so}(p_{2}).
\]
\par
2) We consider $Sp(2n)\cap T^{*}_{r}(n,n).$
Then we can show that
\[
Sp(2n)\cap T^{*}_{r}(n,n)=\left\{ \begin{pmatrix}
 ^{t}C^{-1} & B \\
 O & C
 \end{pmatrix}
\Bigm |\ C\in GL(n,\R), ^{t}\!CB=^{t}\!\!BC\right\},
 \]
which is denoted by $L(2n).$
The condition $^{t}CB=^{t}\!\!\!BC$ means that $^{t}CB$ is a symmetric matrix.
For any symmetric matrix $D\in M_{n}(\R),$ we have $B=^{t}\!\!C ^{-1}D,$
so that we have
\[
L(2n)=\left\{ \begin{pmatrix}
 ^{t}C^{-1} &  ^{t}C^{-1}D\\
 O & C
 \end{pmatrix}
\Bigm |\ C\in GL(n,\R), ^{t}\!\!D=D\right\}.
\] 
It follows that the corresponding Lie algebra is
\[
\mathfrak{l}(2n)=\left\{ \begin{pmatrix}
 -^{t}\!X & Y \\
 O & X
 \end{pmatrix}
\Bigm |\ X\in M_{n}(\R), ^{t}\!Y=Y \right\}.
 \]
Therefore, we can show that $\theta [L(2n)]_{0}(2n)$ is 
\[
\left\{
-\sum_{i=1}^{n}\left(\sum_{k=1}^{n}\frac{\partial \xi_{k}}{\partial y_{i}}(y)x_{k}+
\frac{\partial \eta}{\partial y_{i}}(y)\right)\frac{\partial}{\partial x_{i}}
+\sum_{i=1}^{n}\xi_{i}(y)\frac{\partial }{\partial y_{i}}\ |\ \xi_{i}(y)\in \mathfrak{M}_{n}, \eta(y)\in \mathfrak{M}_{n}^{2}\right\}
\]
and 
\[
\theta [L(2n)](2n)=\left\langle \frac{\partial }{\partial x_{1}},\dots \frac{\partial }{\partial x_{n}}, \frac{\partial }{\partial y_{1}},\dots \frac{\partial }{\partial y_{n}}\right\rangle _{\R}+\theta [L(2n)]_{0}(2n),
\]
where $(x,y)=(x_{1},\dots x_{n},y_{1},\dots ,y_{n})\in (\R^{n}\times \R^{n},0).$

By definition, $\Phi (x,y)\in {\rm Diff}[L(2n)](2n)$ if and only if
$\Phi: (\R^{n}\times\R^{n},0)\lon (\R^{n}\times\R^{n},0)$ is a symplectic diffeomorphism 
of the form $\Phi (x,y)=(\phi_{1}(x,y),\phi_{2}(y))$, where $\phi_{2}\in {\rm Diff}\, (n).$
A symplectic diffeomorphism germ with this property is called a {\it Lagrangian
diffeomorphism germ} in the theory of Lagrangian singularities (cf. \cite[Part III]{Arnold1}).
In this case we can show that $\mathcal{E}_{2n}[L(2n)]=\R.$

\end{Exa}

\section{Infinitesimal structures of geometric equivalence}
{\ \ \ \ }
In this section we now consider $\mathcal{A}[G](n,p)$ for a linear Lie group $G\subset GL(p,\R).$
For a map germ $f:(\R^{n},0)\lon (\R^{p},0),$
we have an $\R$-linear map $$\omega f_{[G]}=\omega f|_{\theta[G]( p)}:\theta[G]( p)\lon \theta (f).$$
We define
\[
T\mathcal{L}[G]_e(f) =\omega f (\theta[G](p )),\ T\mathcal{L}[G](f) =\omega f (\theta[G]_0(p )).
\]
Then $\theta (f)$ is an $\mathcal{E}_p[G]$-module through $f^*_{[G]}=f^*|_{\mathcal{E}_p[G]}:\mathcal{E}_p[G]\lon \mathcal{E}_n$
and $\omega f_{[G]}$ is an $\mathcal{E}_p[G]$-homomorphism over $f^*_{[G]}:\mathcal{E}_p[G]\lon \mathcal{E}_n.$
In this case, $$(\omega f_{[G]}, tf, \theta[G]( p),\theta (n), \theta (f))$$ is the mixed homomorphism over
$f^*_{[G]},$ which is defined in \cite{MatherIII}.
In \S 4 we have considered the case for $G=SL(p,\R)$ or $G=SO(p)$.
Then we have $\mathcal{E}_p[SL(p,\R)]=\R$ and the above mixed homomorphism is not a finite type.
However, we have shown that $\mathcal{E}_p[SO(p)]=\R$ and the above mixed homomorphism is a finite type.

\par
We define
\begin{eqnarray*}
&{}& T\mathcal{L}[G]_e(f)=\omega f_{[G]}(\theta[G](n)), \\
&{}& T\mathcal{L}[G](f)=\omega f_{[G]}(\theta[G]_0(n)), \\
&{}& T\mathcal{R}_e(f)=tf( \theta(n)),\\
&{}& T\mathcal{R}(f)=tf( \mathfrak{M}_n\theta(n)).
\end{eqnarray*}
Then we also define
\begin{eqnarray*}
&{}& T\mathcal{A}[G]_e(f)= T\mathcal{R}_e(f)+T\mathcal{L}[G]_e(f)=tf(\theta (n))+\omega f_{[G]}(\theta[G](n)),\\
&{}& T\mathcal{A}[G](f)=T\mathcal{R}(f)+T\mathcal{L}[G](f)=tf(\mathfrak{M}_n\theta(n))+\omega f_{[G]}(\theta[G]_0(n)).
\end{eqnarray*}
If $G=GL(p,\R),$ then
\[
T\mathcal{A}[GL(p,\R)]_e(f)=T\mathcal{A}_e(f)\ \mbox{and}\ T\mathcal{A}[GL(p,\R)](f)=T\mathcal{A}(f)
\]
in \cite{Wall}.
\begin{Rem}\rm Here, we have a natural question for $\mathcal{A}[G]$: When is $\mathcal{A}[G]$ a geometric subgroup of  $\mathcal{A}$?
It depends on the choice of $G$.
For example if $G=GL(p,\R)$, then $\mathcal{A}[GL(p,\R)]$ is  $\mathcal{A}$ itself, so that
it is a geometric subgroup of $\mathcal{A}.$
There are many examples of $G\subset GL(p,\R)$  such that
$\mathcal{A}[G]$ is not a geometric subgroup of $\mathcal{A}$ (cf. \S 6.2 and \S 6.3).
Since $\mathcal{R}$ is a geometric subgroup of $\mathcal{A},$ the situation depends on $\mathcal{L}[G].$
\end{Rem}
\par
We now focus on $\mathcal{R}\times G$-equivalence.
By definition,  $(\mathcal{R}\times G)(n,p) \subset \mathcal{A}[G](n,p)\subset \mathcal{A}(n,p).$
For a map germ $f:(\R^n,0)\lon (\R^p,0)$,
we set $\mathfrak{g}(f)=\{X.f\ |\ X\in \mathfrak{g}\},$
where $X.f(x)=X.\,^t\!f(x)$ for any $x\in (\R^{n},0).$ 
Then we define the extended tangent spaces and the tangent spaces of the $\mathcal{R}\times G$-orbit  through $f$ by
\[
T(\mathcal{R}\times G)_e(f)=T\mathcal{R}_e(f)+\mathfrak{g}(f),\  
T(\mathcal{R}\times G)(f)=T\mathcal{R}(f)+\mathfrak{g}(f). 
\]
For any $X\in \mathfrak{g},$ we define $\xi_X:(\R^p,0)\lon (\R^p,0)$ by
$\xi_X(y)=X.^{t}y.$ Then $\xi_X$ is a linear mapping, so that it an element of $\theta(p).$
In this sense, we can embed $\mathfrak{g}$ into $\theta (p).$
Therefore we have $\omega f|_\mathfrak{g}:\mathfrak{g}\lon \theta (f),$ which is a $\R$-linear mapping.
This means that $\mathfrak{g}(f)=\omega f(\mathfrak{g}).$
Here $tf: \theta (n)\lon \theta (f)$ is an $\mathcal{E}_n$-homomorphism and $\omega f|_\mathfrak{g}:\mathfrak{g}\lon \mathfrak{g} (f)$
is an $\R$-linear mapping.
Hence, $(\omega f_\mathfrak{g},tf, \mathfrak{g}, \theta (n),\theta (f))$ is the mixed homomorphism of finite type over $f^*|_{\R}=\iota :\R\subset  \mathcal{E}_n.$
Therefore,  $\mathcal{R}\times G$ is a geometric subgroup of $\mathcal{A}$ in the sense of Damon \cite{Damon}.
It sounds a good news, but it is not so good as the following proposition shows.
\begin{Pro}
Let $G\subset GL(p,\R)$ be a Lie subgroup and $f:(\R^n,0)\lon (\R^p,0)$ a map germ.
If $p>1$ and $\dim_\R \theta (f)/T(\mathcal{R}\times G)_e(f)<\infty,$
then $f$ is a submersion germ.
\end{Pro}
\demo Since $T\mathcal{R}_{e}(f)\subset T(\mathcal{R}\times G)_e(f),$ we have
\[
\dim_{\R}\theta (f)/T\mathcal{R}_{e}(f)=\dim_{\R}\theta (f)/T(\mathcal{R}\times G)_e(f)+\dim_{\R}(T(\mathcal{R}\times G)_e(f))/T\mathcal{R}_{e}(f).
\]
Here, $\dim_{\R}(T(\mathcal{R}\times G)_e(f))/T\mathcal{R}_{e}(f)=\dim_{\R}\mathfrak{g}(f)<\infty.$
Therefore, if $\dim_\R \theta (f)/T(\mathcal{R}\times G)_e(f)<\infty,$ then 
$\dim_{\R}\theta (f)/T\mathcal{R}_{e}(f)<\infty.$
\par
On the other hand, it is known (cf. \cite[Proposition 1.11]{MatherR}) that if $p>1$ and $\dim_{\R}\theta (f)/T\mathcal{R}_{e}(f)<\infty,$ then
$f$ is a submersion. 
\enD
For $p=1,$ $GL(1,\R)=\R^{*}=\R\setminus \{0\}$, so that there are only three cases: $G=\{1\}$, $\{\pm 1\},$ or $\R^{*}.$ Therefore, all the cases, classifications by $\mathcal{R}\times G$-equivalence are almost the same 
as the case $G=\{1\}$ (i.e. $\mathcal{R}$-equivalence).
Moreover, if $p>n,$ then there are no finitely determined map germs relative to
$\mathcal{R}\times G.$

\section{Examples of $\mathcal{A}[G]$-equivalence}
\par
{\ \ \ \ \ }
In this section we give some interesting examples of $\mathcal{A}[G]$-equivalence for
$G\subset GL(p,\R).$
We give a survey on the previous results from the view point of our framework on
$\mathcal{A}[G]$-equivalence.

\subsection{Isometric $\mathcal{A}$-equivalence}
{\ \ \ \,}
For $f=(f_{1},f_{2}):(\R^{n},0)\lon (\R^{2},0),$ we have
\[
T\mathcal{L}[SO(2)]_e(f)=\omega f(\theta (SO(2))(2))= \left\langle \frac{\partial }{\partial y_{1}}\circ f,
\frac{\partial }{\partial y_{2}}\circ f\right\rangle_{\R}+\left\langle\left(f_{2}\frac{\partial }{\partial y_{1}}\circ f-f_{1}\frac{\partial }{\partial y_{2}}\circ f\right)\right\rangle _{\R}
\]
and 
\[
T\mathcal{L}[SO(2)](f)=\left\langle\left(f_{2}\frac{\partial }{\partial y_{1}}\circ f-f_{1}\frac{\partial }{\partial y_{2}}\circ f\right)\right\rangle _{\R}.
\]
\par
On the other hand, let $f=(f_1,\dots ,f_p):(\R^n,0)\lon (\R^p,0)$ be a map germ.
By the similar arguments to the above case, we have
\begin{eqnarray*}
T\mathcal{L}[SO( p)]_e(f)\!\!\!&=&\!\!\! \left\langle \left\{\frac{\partial }{\partial y_{i}}\circ f\ \Bigm | 1\leq i\leq p\right\}\right\rangle_{\R} \\
&+&\!\!\! \left\langle\left\{\left(f_{j}\frac{\partial }{\partial y_{i}}\circ f-f_{i}\frac{\partial }{\partial y_{j}}\circ f\right)\ \Bigm |1\leq i < j\leq p\right\}\right\rangle _{\R}
\end{eqnarray*}
and 
\[
T\mathcal{L}[SO( p)](f)=\left\langle\left\{\left(f_{j}\frac{\partial }{\partial y_{i}}\circ f-f_{i}\frac{\partial }{\partial y_{j}}\circ f\right)\ \Bigm |1\leq i < j\leq p\right\}\right\rangle _{\R}.
\]
This is a geometric subgroup of $\mathcal{A}$ in \cite{Damon}. By Proposition 3.5, we have
${\rm Diff}_0[SO( p)](p )=SO(p ).$
Following the classical Euclidean differential geometry, we say that $f,g:(\R^n,0)\lon (\R^p,0)$ are {\it congruent} if there
exists a diffeomorphism germ $\phi:(\R^n,0)\lon (\R^n,0)$ and $A\in SO( p)$ such that
$f\circ \phi (x)=A.g(x)$ for any $x\in (\R^n,0).$
In our terminology $f,g$ are congruent if and only if $f,g$ are $\mathcal{R}\times SO(p )$-equivalent.
By Theorem 4.6, we have the following theorem.
\begin{Th} 
For map germs $f,g:(\R^n,0)\lon (\R^p,0)$, $f, g$ are $\mathcal{A}_{0}[SO( p)]$-equivalent if and only if $f,g$ are congruent.
\end{Th}
This theorem means that the theory of $\mathcal{A}_{0}[SO( p)]$-equivalence among map germs is 
the Euclidean differential geometry on map germs.
Then we have the following corollary of Proposition 5.2 and Theorem 6.1.
\begin{Co} If $p>n$, then there are no map germ $f:(\R^{n},0)\lon (\R^{p},0)$ such that
$\dim_{\R}\theta(f)/T\mathcal{A}[SO(p)]_{e}(f)<\infty.$
\end{Co}
\par
In the case when $n=1,p=2,$ a map germ $f:(\R,0)\lon (\R^2,0)$ is a planer curve germ.
If $f$ is non-singular, we have the curvature function germ $\kappa _f:(\R,0)\lon \R.$
Since the positive or negative sign of the curvature depends on the orientation of the curve, we have $\kappa_{f\circ \phi}(x)={\rm sig}(\phi) \kappa_f(\phi (x))$ for a diffeomorphism germ $\phi :(\R,0)\lon (\R,0),$ where 
\[
{\rm sig}(\phi)=\begin{cases}
+1 &\ {\rm if}\ \dot{\phi}> 0,\\
-1 &\ {\rm if}\ \dot{\phi}< 0 .
\end{cases}
\]
By the classical classification theorem for regular curves in the Euclidean plane $\R^2,$
we have the following proposition.
\begin{Pro} Let $f,g:(\R,0)\lon (\R^2,0)$ be regular map germs. Then $f,g$ are $\mathcal{A}_{0}[SO(2)]$-equivalent if and only if there exists a diffeomorphism germ $\phi:(\R,0)\lon (\R,0)$ such 
that $\kappa _{g}(x)={\rm sig}(\phi) \kappa_f(\phi (x))$ for any $x\in (\R,0).$
\end{Pro} 
For a general map germ $f:(\R,0)\lon (\R^2,0),$ we define its type as follows:
For a function germ $f:(\R,0)\lon (\R,0)$, we say that $f$ has {\it type $A_k$} if $f'(0)=f''(0)=\cdots f^{(k)}(0)=0$
and $f^{(k+1)}(0)\not= 0.$
We also say that $f$ has {\it type $A_{\geq k}$} if $f'(0)=f''(0)=\cdots f^{(k)}(0)=0.$
We have the following lemma \cite[Theorem 3.3]{Bru-Gib}.
\begin{Lem}
Let $f:(\R,0)\lon (\R,0)$ be a function germ of type $A_k$. Then there exists a diffeomorphism germ $\phi :(\R,0)\lon (\R,0)$
such that $f\circ \phi (x)=\pm x^{k+1},$ where $+$ or $-$ according as $f^{(k+1)}(0)$ is positive or negative.
\end{Lem}
We now consider a map germ $f:(\R,0)\lon (\R^2,0)$ with $f(x)=(f_1(x),f_2(x)).$
We say that $f$ has {\it type $A_k$} if one of $f_1$ or $f_2$ has type $A_k$ and another has type $A_{\geq k}$. 
Then we have the following proposition.
\begin{Pro} Suppose $f:(\R,0)\lon (\R^2,0)$ has type $A_k$. Then there exists a function germ $h:(\R,0)\lon \R$ such that $f$ is $\mathcal{A}_{0}[SO(2)]$-equivalent to
$(\pm x^{k+1},x^{k+1}h(x)).$ 
\end{Pro} 
\demo
For $f=(f_1,f_2),$ suppose $f_1$ has type $A_k$ and $f_2$ has type $A_{\geq k}.$
By Lemma 6.4, there exists a diffeomorphism germ $\phi :(\R,0)\lon (\R,0)$ such that $f_1\circ \phi (x)=\pm x^{k+1}.$
By the Hadamard lemma, there exists $h$ such that $f_2\circ \phi (x)=x^{k+1}h(x).$
Instead, if $f_2$ has type $A_k$ and $f_1$ has type $A_{\geq k},$ 
there exists a diffeomorphism germ $\phi :(\R,0)\lon (\R,0)$ and a function germ 
$g:(\R,0)\lon \R$ such that $f\circ \phi (x)=(x^{k+1}g(x),\pm x^{k+1})$
by the similar reason to above.
Then we have
\[
\begin{pmatrix}
0 & 1 \\
-1 & 0
\end{pmatrix}
\begin{pmatrix}
x^{k+1}g(x) \\
\pm x^{k+1}
\end{pmatrix}
=
\begin{pmatrix}
\pm x^{k+1} \\
x^{k+1}(-g(x))
\end{pmatrix}.
\]
If we put $h(x)=-g(x),$ then
$f$ is $\mathcal{A}[SO(2)]$-equivalent to $(\pm x^{k+1},x^{k+1}h(x)).$
Since $\begin{vmatrix} 0 & 1 \\ -1 & 0\end{vmatrix}=1,$ it is actually $\mathcal{A}_{0}[SO(2)]$-equivalence.
\enD
\par
\begin{Rem}\rm
In \cite{Bru-Gaff} a classification of $\mathcal{A}$ simple map germs $(\mathbb{C},0)\lon (\mathbb{C}^2,0)$ has been given.
Here $\mathcal{A}=\mathcal{A}[GL(2,\mathbb{C})].$ They have shown that $f$ is $\mathcal{A}$-equivalent to
$(x^{k+1},x^{k+1}h(x))$ and $f$ is not $\mathcal{A}$ simple if $k>3$ or $k=3$ and $h$ has type $A_{\geq 4}.$ 
However, for the classification by $\mathcal{A}[SO(2)]$-equivalence, there might be no $\mathcal{A}[SO(2)]$ simple germs.
\end{Rem}
\par
We consider a map germ $f:(\R,0)\lon (\R^2,0)$ with type $A_k.$
By Proposition 6.5, we may assume that $f(x)=(\pm x^{k+1},x^{k+1}h(x)).$
In this case we have
\[
\dot{f}(x)=(\pm (k+1)x^k, x^k((k+1)h(x)+x\dot{h}(x))),
\]
so that the singular point of $f$ is the origin.
We define
\[
\bm{\mu}(x)=\frac{1}{\sqrt{(k+1)^2+((k+1)h(x)+x\dot{h}(x))^2}}(\pm (k+1), (k+1)h(x)+x\dot{h}(x)).
\]
Then $\bm{\mu}(x)$ is a unit vector tangent to $f(\R)$ at $x\not= 0$ (i.e. a regular point of $f$).
We also define $\bm{\nu}(x)=J\bm{\mu}(x)$, where $J=\begin{pmatrix} 0 & -1 \\ 1 & 0 \end{pmatrix}.$
It follows that $\{ \bm{\nu}(x), \bm{\mu}(x)\}$ is an orthonormal frame along $f.$
Moreover, we have a map germ $(f,\bm{\nu}):(\R,0)\lon \R^2\times S^1$ with $\dot{f}(x)\cdot \bm{\nu}(x)=0,$ where
$\bm{a}\cdot \bm{b}$ is the canonical scaler product of $\R^2.$
This means that $f$ is a frontal in the sense of \cite{Fuku-Taka}. If we define $\ell_f(x)=\dot{\bm{\nu}}(x)\cdot \bm{\mu}(x)$
and $\beta_f(x)=\dot{f}(x)\cdot\bm{\mu}(x),$ we have the following Frenet-type formulae \cite{Fuku-Taka}:
\[
\begin{pmatrix} \dot{\bm{\nu}}(x) \\ \dot{\bm{\mu}}(x) \end{pmatrix}=\begin{pmatrix} 0 & \ell _f(x) \\ -\ell _f(x) & 0 \end{pmatrix}\begin{pmatrix} \bm{\nu}(x)\\ \bm{\mu}(x) \end{pmatrix},\ \dot{f}(x)=\beta_f(x)\bm{\mu}(x).
\]
The following uniqueness theorem was shown in \cite{Fuku-Taka}:
\begin{Th}
Let $f,g:(\R,0)\lon (\R^2,0)$ be frontal germs. Then $f,g$ are congruent if there exists a diffeomorphism germ
$\phi:(\R,0)\lon (\R,0)$  such that $\dot{\phi}(x)>0$ and
$$
\dot{\phi}(x)(\ell _f\circ\phi (x), \beta_f\circ \phi (x))=(\ell _g(x), \beta_g(x))
$$ for any $x\in (\R,0).$ 
\end{Th}
As a special case, we have a classification theorem on map germ $f:(\R,0)\lon (\R^2,0)$ with type $A_k.$
For $f(x)=(\pm x^{k+1},x^{k+1}h(x)),$ we can show that
\begin{eqnarray*}
\ell_f(x)&=&\frac{\pm (k+1)((k+2)\dot{h}(x)+x\ddot{h}(x))}{\sqrt{((k+1)^2+((k+1)h(x)+x\dot{h}(x))^2)^3}}, \\
\beta_f(x)&=& x^k\sqrt{(k+1)^2+((k+1)h(x)+x\dot{h}(x))^2}.
\end{eqnarray*}
Therefore, the basic invariant $(\ell_f,\beta_f)$ depends on $h(x).$
\par
On the other hand, we consider $n=2$ and $p=3.$
For a regular surface, we have the Monge normal form.
By the classification theorem for quadratic forms, we have the following proposition.
\begin{Pro}
Let $f:(\R^{2},0)\lon (\R^{3},0)$ be an immersion germ. Then $f$ is
$\mathcal{R}\times SO(3)$-equivalent {\rm (}i.e. $\mathcal{A}_{0}[SO(3)]$-equivalent\/{\rm )} to the following germ:
\[
g(x_{1},x_{2})=(x_{1},x_{2},\lambda_{1}x_{1}^{2}+\lambda_{2}x_{2}^{2}+a_{30}x_{1}^{3}+a_{21}x_{1}^{2}x_{2}+a_{12}x_{1}x_{2}^{2}+a_{03}x_{2}^{3}+O(4)).
\]
Here, $\lambda _{1},\lambda_{2}$ are the principal curvatures of $f$ at the origin.
\end{Pro}
The above map germ $g$ is called a {\it Monge normal form}.
Recently, $\mathcal{R}\times SO(3)$-equivalence has been used for the study of differential geometry of singular surfaces in $\R^3$ (cf. \cite{F-H, F-T2, M-S, M-S2,West}).

\subsection{Volume preserving $\mathcal{A}$-equivalence}
{\ \ \ \ \ }
We now consider the case when $G=SL(p,\R).$
In this case, $\mathcal{A}[SL(p,\R)]$-equivalence is volume preserving $\mathcal{A}$-equivalence on
the target space.
By Example 4.7, $\theta [SL(p,\R)](p)$ is not a finitely generated
$(\mathcal{E}_{p}[SL(p,\R)]=\R)$-module, so that $\mathcal{L}[SL(p,\R)]$ and $\mathcal{A}[SL(p,\R)]$ are not geometric subgroups of $\mathcal{A}$ in the sense of Damon \cite{Damon}.
Therefore, the usual techniques of the singularity theory cannot work properly.
However, as Martinet (cf. \cite[page 50]{martinet}) pointed out, the group ${\rm Diff}[SL(p,\R)](p)$ is big enough that there is still some hope of finding a reasonable classification theorem by volume preserving $\mathcal{A}$-equivalence.
Actually, Domitrz and Rieger investigated this equivalence in \cite{D-R}. 
In their paper $\mathcal{A}[SL(p,\R)]$ is written by $\mathcal{A}_{\Omega_{p}}.$
They called the geometry associated with ${\rm Diff}[SL(p,\R)](p)$ a {\it unimodular geometry}.
For convenience, they adopted $\mathcal{A}[SL(p,\mathbb{C})]$ instead of $\mathcal{A}[SL(p,\R)].$
One of their classifications is as follows.
\begin{Pro}[\cite{D-R}] 
Any $\mathcal{A}[SL(2,\mathbb{C})]$-simple map-germ $f:(\mathbb{C}^{n},0)\lon (\mathbb{C}^{2},0)$,
$n\geq 2,$
is $\mathcal{A}[SL(2,\mathbb{C})]$-equivalent to one of the following list of germs\/{\rm :}
$$
(x_{1},x_{2});\ (x_{1},x_{2}^{2}+Q);\ (x_{1},x_{1}x_{2}+x_{2}^{3}+Q);\ (x_{1},x_{2}^{3}+x_{1}^{k}x_{2}+Q),k>1;
\ (x_{1}x_{1}x_{2}+x_{2}^{4}+Q),
$$
where $Q=\sum_{i=3}^{n}x_{i}^{2}$ for $n>2$ and $Q=0$ for $n=2.$
\end{Pro}
Since $SL(2,\mathbb{C})=Sp(2,\mathbb{C}),$ the above list also gives a classification of
$\mathcal{A}[Sp(2,\mathbb{C})]$-simple map germs $(\mathbb{C}^{n},0)\lon (\mathbb{C}^{2},0)$,
$n\geq 2.$ For $n=1,$ Ishikawa and Janeczko classified $\mathcal{A}[Sp(2,\mathbb{C})]$-simple map germs
$(\mathbb{C},0)\lon (\mathbb{C}^{2},0)$ in \cite{Ishi-Jane}.
However, classifications by $\mathcal{A}[Sp(2p)]$-equivalence for general $p$ is much complicated.
\par
Moreover,
$(\mathcal{R}\times SL(p,\R))(n,p)$ is a proper subgroup of $\mathcal{A}[SL(p,\R)](n,p).$
The notion of $\mathcal{R}\times SL(p,\R)$-equivalence is known to be
{\it equi-affine congruence.}
We say that $f,g:(\R^{n},0)\lon (\R^{p},0)$ are {\it equi-affine congruent} if
there exist $\phi \in {\rm Diff}\,(n)$ and $A\in SL(p,\R)$ such 
that $f\circ \phi (x)=A.g(x)$ for any $x\in (\R^{n},0).$
For a regular curve $f:(\R,0)\lon (\R^{2},0)$
without inflection points, an {\it equi-affine curvature} of $f$ is defined and it is denoted
by $\kappa^{e}_{f}$ (cf. \cite{Nomizu-Sasaki}). The following uniqueness theorem is known.
\begin{Pro} Let $f,g:(\R,0)\lon (\R^2,0)$ be regular map germs without inflection points. Then $f,g$ are $\mathcal{R}\times SL(2,\R)$-equivalent {\rm (}i.e. equi-affine congruent\/{\rm )} if and only if there exists a diffeomorphism germ $\phi:(\R,0)\lon (\R,0)$ such 
that $\kappa ^{e}_{g}(x)={\rm sig}(\phi) \kappa^{e}_f(\phi (x))$ for any $x\in (\R,0).$
\end{Pro}
The equi-affine geometry for singular curves is also an interesting subject.
Moreover, the surface theory (i.e. $n=2,p=3$) is also quite interesting.
These are our future assignments.
As a consequence, there is a big gap between the unimodular geometry and the equi-affine geometry.
This is completely different from the case when $G=SO(p).$

\subsection{Bi-$\mathcal{A}$-equivalence}
{\ \ \ \ \ }
We consider a map germ $f=(f_1,f_2):(\R^n,0)\lon (\R^{p_{1}}\times\R^{p_{2}},0)$
which is considered to be a divergent diagram of map germs $(\R^{p_{1}},0)\stackrel{f_{1}}{\longleftarrow} (\R^{n},0)\stackrel{f_{2}}{\longrightarrow} (\R^{p_{2}},0).$
The notion of bi-$\mathcal{A}$-equivalence among map germs of the form $f=(f_1,f_2):(\R^n,0)\lon (\R^{p_1}\times \R^{p_2},0)$
was introduced in \cite{CPS,Dufour1,Ped}. 
We say that $f=(f_1,f_2)$ and $g=(g_1,g_2)$ are {\it bi-$\mathcal{A}$-equivalent} if there 
exist diffeomorphism germs $\phi :(\R^n,0)\lon (\R^n,0)$ and $\psi_i:(\R^p,0)\lon (\R^{p_i},0)$, $(i=1,2)$, such that
$f_i\circ \phi=\psi_i\circ g_i.$
We consider the Lie group $D^{*}(p_{1},p_{2})\subset GL(p,\R),$ where
$p=p_{1}+p_{2}.$ 
Then $f=(f_1,f_2)$ and $g=(g_1,g_2)$ are $\mathcal{A}[D^{*}(p_{1},p_{2})]$-equivalent if and only if
these are bi-$\mathcal{A}$-equivalent.

In this case, we have
$\theta [D^{*}(p_{1},p_{2})](p )=\theta (p_{1})\oplus\theta (p_{2}).$ 
Then
\[
T\mathcal{A}[D^{*}(p_{1},p_{2})]_e(f)=tf ( \theta (n))+\omega f(\theta(p_{1})\oplus \theta (p_{2})),
\]
\[
T\mathcal{A}[D^{*}(p_{1},p_{2})](f)=tf ( \mathfrak{M}_n\theta (n))+\omega f(\mathfrak{M}_{p_{1}}\theta(p_{1})\oplus \mathfrak{M}_{p_{2}}\theta (p_{2})).
\]
This is not a geometric subgroup of $\mathcal{A}$ in the sense of Damon \cite{Damon}.
In particular, if we consider the case $p=2,$ the bi-$\mathcal{A}$-stable map germs
were classified by Dufour \cite{Dufour1}. Moreover, a formal classification for formal finite bi-$\mathcal{A}$-codimensional map germs was given by Mancini, Ruas and Texieira \cite{Mancini}.

\subsection{Strict bi-$\mathcal{A}$-equivalence }
{\ \ \ \ \ }
Here we also consider divergent diagrams
$(\R^{p_{1}},0)\stackrel{f_{1}}{\longleftarrow} (\R^{n},0)\stackrel{f_{2}}{\longrightarrow} (\R^{p_{2}},0).$
We say that $f=(f_1,f_2)$ and $g=(g_1,g_2)$ are {\it strictly bi-$\mathcal{A}$-equivalent} if there 
exist diffeomorphism germs $\phi :(\R^n,0)\lon (\R^n,0)$ and $\psi_2:(\R^p,0)\lon (\R^{p_2},0)$ such that
$f_{1}\circ \phi=g_{1}$ and $f_2\circ \phi=\psi_2\circ g_2.$
Then we consider the following Lie group:
$$\{I_{p_1}\}\oplus GL(p_2,\R)=\left\{\begin{pmatrix} I_{p_1} & 0  \\
0 & A 
\end{pmatrix}\Bigm |A\in GL(p_2,\R) \right\}
$$
In this case, we have
$\theta [\{I_{p_1}\}\oplus GL(p_2,\R)](p)=\{0\}\oplus \theta (p_2),$ so that
\[
T\mathcal{A}[\{I_{p_1}\}\oplus GL(p_2,\R)]_e(f)=tf ( \theta (n))+\omega f(\{0\}\oplus \theta (p_2)),
\]
\[
T\mathcal{A}[\{I_{p_1}\}\oplus GL(p_2,\R)](f)=tf ( \mathfrak{M}_n\theta (n))+\omega f(\{0\}\oplus \mathfrak{M}_{p_2}\theta (p_2)).
\]
This is a geometric subgroup of $\mathcal{A}$ in the sense of Damon \cite{Damon}.
By definition, $f=(f_1,f_2),g=(g_1,g_2):(\R^n,0)\lon (\R^{p_1}\times \R^{p_2},0)$ are $\mathcal{A}[\{I_{p_1}\}\oplus GL(p_2,\R))]$-equivalent if and only if these germs are strictly bi-$\mathcal{A}$-equivalent.
However, this equivalence is too strong.
If $f=(f_1,f_2),g=(g_1,g_2)$ are strictly bi-$\mathcal{A}$-equivalent, then
$f_{1},g_{1}$ are $\mathcal{R}$-equivalent.
By the same reason as Proposition 5.2 (cf. \cite[Proposition 1.11]{MatherR}),
$f_{1}$ is a submersion for $p_{1}>1$ and
$\dim_{\R}\theta (f)/T\mathcal{A}[\{I_{p_1}\}\oplus GL(p_2,\R)]_e(f)<\infty.$
We remark that even if $p_{1}=1,$ 
functional moduli appear for very low dimensions (i.e. $p_{2}=2$, \cite{Dufour}).
In order to avoid functional moduli, we consider
another Lie group defined by
\[
(1^+,GL(p_{2},\R))=\left\{\begin{pmatrix} 1 & \bm{b}  \\
0 & A 
\end{pmatrix}\Bigm |\bm{b}\in \R^{p_{2}},\ A\in GL(p_{2},\R) \right\}\subset GL(p,\R).
\]
Then $f,g$ are $\mathcal{A}[(1^+,GL(p_2,\R))]$-equivalent if and only if 
there exist  diffeomorphism germs $\phi :(\R^n,0)\lon (\R^n,0)$,
$\psi :(\R^{p_2},0)\lon (\R^{p_2},0)$ and a function germ $\alpha :(\R^{p_2},0)\lon (\R,0)$
such that $f_1(x)+\alpha (f_2(x))=g_1\circ \phi (x)$ and 
$
\psi\circ f_2=g_2\circ\phi .
$
In this case, 
the corresponding Lie algebra is
\[
(0^{+},M_{p_{2}}(\R))=\left\{\begin{pmatrix} 0 & \bm{b}  \\
0 & X 
\end{pmatrix}\Bigm |\bm{b}\in \R^{p_{2}},\ X\in M_{p_{2}}(\R) \right\}\subset M_{p}(\R), 
\]
so that we have
$$
\theta [(1^+,GL(p_{2},\R))](p)=\left\langle \frac{\partial}{\partial x},\frac{\partial}{\partial y_{1}},
\dots \frac{\partial}{\partial y_{p_{2}}}\right\rangle _{\mathcal{E}_{p_{2}}},
$$
where $(x,y_{1},\dots y_{p_{2}})\in \R\times\R^{p_{2}}=\R^{p}.$
Therefore, we have
\[
T\mathcal{A}[(1^+,GL(p_{2},\R))]_e(f)=tf ( \theta (n))+\omega f\left(\left\langle \frac{\partial}{\partial x},\frac{\partial}{\partial y_{1}},
\dots \frac{\partial}{\partial y_{p_{2}}}\right\rangle _{\mathcal{E}_{p_{2}}}\right),
\]
\[
T\mathcal{A}[(1^+,GL(p_{2},\R))](f)=tf ( \mathfrak{M}_{n}\theta (n))+\omega f\left(\left\langle \frac{\partial}{\partial x},\frac{\partial}{\partial y_{1}},
\dots \frac{\partial}{\partial y_{p_{2}}}\right\rangle _{\mathfrak{M}_{p_{2}}}\right).
\]
This is a geometric subgroup of $\mathcal{A}$ in the sense of Damon \cite{Damon}.
In \cite{Dufour} a generic classification of $f=(f_{1},f_{2}):(\R^2,0)\lon (\R^3,0)=(\R\times \R^{2},0)$ with respect to $\mathcal{A}[\{1\}\oplus GL(2,\R)]$-equivalence was given. One of the normal form is
\[
f_{1}(x_{1},x_{2})=\pm x_{1}+\alpha\circ f_{2}(x_{1},x_{2}),\ f_{2}(x_{1},x_{2})=(x_{1}^{3}+x_{1}x_{2},x_{2}).
\]
Here $\alpha :(\R^{2},0)\lon \R$ is the functional modulus.
By definition, it is $\mathcal{A}[(1^+,GL(p_2,\R))]$-equivalent to
$(\pm x_{1},x_{1}^{3}+x_{1}x_{2},x_{2}).$
Therefore, $\mathcal{A}[(1^+,GL(p_2,\R))]$-equivalence is a strict bi-$\mathcal{A}$-equivalence modulo
functional uni-moduli.
It is known that the functional moduli play an important role in the
geometry of webs.

\subsection{Projections of map germs}
\subsubsection{General projection}
{\ \ \ \ \ }
We consider map germ $f=(f_{1},f_{2}):(\R^{n},0)\lon (\R^{p}=\R^{p_{1}}\times\R^{p_{2}},0)$
and the canonical projection $\pi_{2}:(\R^{p_{1}}\times\R^{p_{2}},0)\lon (\R^{p_{1}},0).$
We say that $f,g:(\R^{n},0)\lon (\R^{p_{1}}\times\R^{p_{2}},0)$ are {\it projection $\mathcal{A}$-equivalent}
with respect to $\pi_{2}$ if there exist $\phi\in {\rm Diff}\,(n)$ and
$\Psi\in {\rm Diff}\, (p)$ of the form $\Psi (x,y)=(\psi_{1}(x,y),\psi_{2}(y))$ such that
$\Psi\circ f=g\circ \phi.$
Then we consider the Lie group $T^{*}_{r}(p_{1},p_{2})\subset GL(p,\R).$
We can show that $f,g:(\R^{n},0)\lon (\R^{p_{1}}\times\R^{p_{2}},0)$ are projection $\mathcal{A}$-equivalent
with respect to $\pi_{2}$ if and only if these are $\mathcal{A}[T^{*}_{r}(p_{1},p_{2})]$-equivalent.
By Example 4.10, 2), we have
\[
T\mathcal{A}[T^{*}_{r}(p_{1},p_{2})]_{e}(f)=tf(\theta (n))+\omega f(\theta (\pi_{p_{1}})\oplus \theta (p_{2})),
\]
\[
T\mathcal{A}[T^{*}_{r}(p_{1},p_{2})](f)=tf(\mathfrak{M}_{n}\theta (n))+\omega f(\mathfrak{M}_{p}\theta (\pi_{p_{1}})\oplus \mathfrak{M}_{p_{2}}\theta (p_{2})).
\]
This is a geometric subgroup of $\mathcal{A}$ in the sense of Damon \cite{Damon}.
We expect that there might be interesting properties on this equivalence.
For example, one of the results of Romero-Fuster, Mancini and Soares-Ruas \cite[Lemma 2]{FMR} is
interpreted by using this equivalence as follows:
\begin{Pro} Suppose that $f=(f_1,f_2),g=(g_1,g_2): (\R^{n},0)\lon (\R^{p_{1}}\times\R^{p_{2}},0)$
are immersion germs.
Then $f_{2},g_{2}$ are $\mathcal{A}$-equivalent if and only if
$f,g$ are $\mathcal{A}[T^{*}_{r}(p_{1},p_{2})]$-equivalent.
\end{Pro}
Since $\mathcal{A}$-equivalence among projections of surfaces (i.e. $n=2,p=3$)
is a useful tool for the study of differential geometry of surface from the view point of contact with lines or planes,
there might be interesting applications of $\mathcal{A}[T^{*}_{r}(p_{1},p_{2})]$-equivalence for general singular surfaces.
Moreover, Mancini and Soares-Ruas \cite{M-R} investigated the following equivalence relation:
Two map germs $f=(f_1,h), g=(g_2,h): (\R^{n},0)\lon (\R^{p_{1}}\times\R^{p_{2}},0)$
are {\it $h$-equivalent} if there exists $\phi\in {\rm Diff}\,(n)$, $\Psi\in {\rm Diff}\, (p_1,p_2)$
of the form $\Psi (y,z)=(\psi_1(y,z),\psi_2(z))$ for $(y,z)\in \R^{p_1}\times\R^{p_2}$ such 
that $\Psi\circ f=g\circ \phi.$
By definition, $f=(f_1,h), g=(g_2,h)$ are $h$-equivalent if and only if
these are $\mathcal{A}[T^{*}_{r}(p_{1},p_{2})]$-equivalent.
\subsubsection{Lagrangian equivalence}
{\ \ \ \,\,}
We now consider a symplectic manifold $\R^{n}\times \R^{n}$ with the canonical
symplectic structure $\omega=\sum_{i=1}^{n} dx_{i}\wedge dy_{i}$
where $(x,y)=(x_{1},\dots ,x_{n},y_{1},\dots ,y_{n}).$
We say that two map germs $f,g:(\R^{m},0)\lon (\R^{n}\times \R^{n},0)$
are {\it Lagrangian equivalent} if there exist $\phi\in {\rm Diff}\,(m)$ 
and a symplectic diffeomorphism $\Psi:(\R^{n}\times \R^{n},0)\lon (\R^{n}\times \R^{n},0)$
of the form $\Psi (x,y)=(\psi_{1}(x,y),\psi_{2}(y))$ such that
$\Psi\circ f=g\circ \phi.$
By Example 4.11, $f,g$ are Lagrangian equivalent if and only if
these are $\mathcal{A}[L(2n)]$-equivalent.
In this case, we have
\[
T\mathcal{A}[L(2n)]_{e}(f)=tf(\theta (m))+\omega f(\theta [L(2n)]),
\]
\[
T\mathcal{A}[L(2n)](f)=tf(\mathfrak{M}_{m}\theta (m))+\omega f(\theta[L(2n)]_{0}),
\]
where
$\theta [L(2n)]_{0}$ and $\theta [L(2n)]$ are given in Example 4.11.
\par
We call $f:(\R^{n},0)\lon (\R^{n}\times \R^{n},0)$ an {\it isotropic map germ}
if $f^{*}\omega =0.$ Moreover, $f:(\R^{n},0)\lon (\R^{n}\times \R^{n},0)$ is said to be a {\it
Lagrangian immersion germ} if $f$ is isotropic and an immersion germ.
For Lagrangian immersion germs, there is a theory of generating families.
A generic classification of Lagrangian immersion germs by using generating families is
well-known (cf. \cite[Part III]{Arnold1}).
However, we do not know classifications of general map germs $f:(\R^{m},0)\lon (\R^{n}\times \R^{n},0)$
by $\mathcal{A}[L(2n)]$-equivalence
so far.
\section{
$A[G]$-geometry versus classical $G$-geometry}
The basic tools of the (local) theory of $\mathcal{A}$-equivalence are 
finite determinacy of map germs and the versality theorem for unfoldings of map germs
which are characterized by the algebraic structure of the formal tangent space
of the $\mathcal{A}$-equivalence class.
In \cite{Damon} Damon gave a very wide class of subgroups of $\mathcal{A}(n,p)$
for which the finite determinacy theorem and the versality theorem hold.
Those subgroups are called {\it geometric subgroups of} $\mathcal{A}.$
However, as we already mentioned in this paper, $\mathcal{A}[G](n,p)$ for
some $G\subset GL(p,\R)$ are not geometric subgroups of $\mathcal{A}$.
\par
On the other hand, we defined $\mathcal{R}\times G$-equivalence among
map germs (cf. \S 2). From the view point of extrinsic differential geometry,
we say that map germs $f,g:(\R^n,0)\lon (\R^p,0)$ are
{\it $G$-congruent} if they are $\mathcal{R}\times G$-equivalent.
We call the corresponding geometry a {\it classical $G$-{\rm (}differential\/{\rm )}geometry}.
Since $(\mathcal{R}\times G)(n,p)\subset \mathcal{A}[G](n,p),$ 
a geometry corresponding to $\mathcal{A}[G]$-equivalence is called an {\it $\mathcal{A}[G]$-geometry}
(or, a {\it differential $G$-geometry}\/) which is categorically wider 
than the classical $G$-geometry.
The classical Euclidean differential geometry is a geometry which investigates invariant quantities
 and properties of immersions $f:(\R^n,0)\lon (\R^p,0)$ under
$SO(p)$-congruence,  
so that it is the classical $SO(p)$-geometry.
It is known that the curvatures of regular curves in $\R^2$ are complete invariants
in the Euclidian differential geometry, which is also
the functional moduli with respect to $SO(2)$-congruence (cf. Proposition 6.3).
We have shown in Theorem 6.1 that the $\mathcal{A}_0[SO(p)]$-geometry 
and the classical $SO(p)$-geometry are the same.
From this point of view, the $\mathcal{A}[GL(p,\R)]$-geometry
is the $\mathcal{A}$-geometry, which should be called a {\it local differential topology} for
map germs.
The classical $GL(p,\R)$-geometry is usually called a {\it full affine geometry}.
Moreover, the $\mathcal{A}[SL(p,\R)]$-geometry is called the {\it unimodular geometry}
in \cite{D-R}.
Since $SL(p,\R)$ is big enough as a Lie subgroup of $GL(p,\R),$
there are not so much difference between the $\mathcal{A}[SL(p,\R)]$-geometry
and the $\mathcal{A}[GL(p,\R)]$-geometry (cf. Proposition 6.9).
The classical $SL(p,\R)$-geometry is known to be the equi-affine geometry (cf. \cite{Nomizu-Sasaki}).
In this case, the equi-affine curvatures for regular curves without inflections in $\R^2$
are also the complete invariants and the functional moduli with respect to equi-affine congruence (i.e. $SL(p,\R)$-congruence).
\par
Since the formal tangent space of an $\mathcal{A}[G]$-equivalence class gives several information,
the calculation of tangent space might be the first step for the investigation of the $\mathcal{A}[G]$-geometry of map germs.
For example, we have the following theorem which is a generalization of Theorem 6.1.
\begin{Th} Let $G\subset GL(p,\R)$ be a connected linear Lie group.
Then the following conditions are equivalent\/{\rm :}
\par\noindent
{\rm (1)} $\mathfrak{g}\cong \theta [G]_0(p)$ as $\R$-vector spaces,
\par\noindent
{\rm (2)} For any $\eta (y)=\sum_{i=1}^p\eta_i(y)(\partial /\partial y_i)\in \theta [G]_0(p),$
$\eta_i (y)$ $(i=1,\dots ,p)$ are linear functions,
\par\noindent
{\rm (3)} ${\rm Diff}_0[G](p)=G,$
\end{Th}
\demo
We consider an $\R$-linear mapping $\iota :M_p(\R)\lon \mathfrak{M}_p\theta (p)$ defined by
$\iota (e_{ij})=y_i(\partial/\partial y_j),$ where $\{e_{ij}\ |\ i,j=i,\dots p\}$ is the
canonical basis of $M_p(\R).$
Since $\{y_i(\partial /\partial y_j)\ |\ i,j=1,\dots ,p\}$ are linearly independent,
$\iota$ is a monomorphism.
Therefore, $\mathfrak{g}\cong \iota(\mathfrak{g})\subset \theta [G]_0(p)$, generally.
Let $\{\delta _1,\dots \delta_r\}$ be a basis of $\mathfrak{g}.$
Then $\delta _i$ is denoted by a linear combination of $e_{ij}.$
If $\mathfrak{g}\cong \theta [G]_0(p),$ then any $\eta (y)=\sum_{k=1}^p\eta _k(y)(\partial /\partial y_k)\in \theta [G]_0(p)=\iota (\mathfrak{g})$ is denoted by a linear 
combination of $y_i(\partial /\partial y_j),$ so that $\eta _k(y)$ are linear functions.
Suppose that $\eta_i (y)$ $(i=1,\dots ,p)$ are linear functions, for any $\eta (y)=\sum_{i=1}^p\eta_i(y)(\partial /\partial y_i)\in \theta [G]_0(p).$
Then $\eta _i(y)=\sum_{k=1}^p a_{ik}y_k.$ By definition $(a_{ik})=((\partial \eta _i/\partial y_k)(y))
\in \mathfrak{g}$ for any $y\in (\R^p,0).$
This means that $\eta (y)=\sum_{i,k=1}^p a_{ik}y_k(\partial/\partial y_i)=\iota (\sum_{i,k=1}^p a_{ik}e_{ki})\in \iota (\mathfrak{g}).$
Thus we have $\iota (\mathfrak{g})=\theta [G]_0(p).$
We have shown that conditions (1) and (2) are equivalent.
\par
We assume that conditions (1) and (2).
Then we can show that condition (3) holds exactly the same method of the proof for Theorem 4.6.
If we consider the formal tangent space of ${\rm Diff}_0[G](p)=G,$ we can easily show that condition (3) implies condition (1).
\enD
\begin{Co} Suppose that one of the conditions in Theorem 7.1 holds. Then we have the following\/{\rm :}
\par\noindent
{\rm (1)} $T\mathcal{A}[G](f)=T(\mathcal{R}\times G)(f)$ for any $f:(\R^n,0)\lon (\R^p,0),$
\par\noindent
{\rm (2)} $f,g :(\R^n,0)\lon (\R^p,0)$ are $\mathcal{A}_0[G]$-equivalent if and only if 
these are $G$-congruent.
\end{Co}
\demo
If $\mathfrak{g}\cong \theta [G]_0(p)$ as $\R$-vector spaces, then $\omega f(\theta [G]_0(p))=\mathfrak{g}(f)$ for any $f:(\R^n,0)\lon (\R^p),$ so that assertion (1) holds.
If $f,g$ are $\mathcal{A}_0[G]$-equivalent, then there exists $(\phi,\psi)\in {\rm Diff}\, (n)\times {\rm Diff}_0[G](p)$ such that $\psi \circ f=g\circ \phi.$ 
Since ${\rm Diff}_0[G](p)=G,$ $f,g$ are $\mathcal{R}\times G$-equivalent (i.e. $G$-congruent).
The converse assertion holds by definition.
\enD
\par\noindent
Corollary 7.2 shows that if one of the conditions in Theorem 7.1 holds,
then the $A[G]$-geometry and the classical $G$-geometry are the same.
For a Lie subgroup $H<G<GL(p,\R),$ we have $\theta [H]_0(p)\subset \theta [G]_0(p).$
If condition (2) in Theorem 7.1 holds for $G,$ then it also holds for $H.$
If $G=SO(p_1,p_2),$ then we can show that condition (2) in Theorem 7.1 holds, where
we do not give the proof of this fact except in the case when $p_1=0, p_2=p$ (cf. Example 4.4).
For $G=SO(4),$ we consider $H=\{1\}\oplus SO(3)$.
Then $f=(f_1,f_2), g=(g_1,g_2) : (\R^2,0)\lon (\R\times \R^3,0)$ are
$\mathcal{A}_0[\{1\}\oplus SO(3)]$-equivalent if and only if 
these are $\mathcal{R}\times (\{1\}\oplus SO(3))$-equivalent,
which also means that there exists $\phi \in {\rm Diff}\, (2)$ and $A\in SO(3)$ such that
$f_1\circ \phi =g_1$ and $f_2\circ \phi=A.g_2.$
Since $f_1\circ \phi=g_1,$ we have $\phi(g_1^{-1}(c))=f_1^{-1}(c)$ for any $c\in (\R,0).$
Therefore, if $f_2,g_2$ are immersive, then the classical $(\{1\}\oplus SO(3))$-geometry
is the classical Euclidean geometry among regular surfaces with singular foliations.
On the other hand, the local differential topology among surfaces with singular foliations is 
the $\mathcal{A}[\{1\}\oplus GL(3,\R)]$-geometry which is different from
the full affine geometry among surfaces with singular foliations
(i.e. the classical $(\{1\}\oplus GL(3,\R))$-geometry).

\par
Moreover, we define an $\R$-vector space
\[
\mathcal{M}(\mathcal{A}[G];\mathcal{R}\times G)(f)=\frac{T\mathcal{A}[G](f)}{T(\mathcal{R}\times G)(f)},
\]
which is called a {\it relative infinitesimal moduli space of} $f$ {\it with respect to
$\mathcal{A}[G]$ and $\mathcal{R}\times G$.}
By the calculations of tangent spaces of $\mathcal{A}[G]$-equivalence and $\mathcal{R}\times G$-equivalence, the moduli space clarifies the difference between
the $\mathcal{A}[G]$-geometry and the classical $G$-geometry.
By Theorems 4.6 and 6.1, $\mathcal{M}(\mathcal{A}[SO(p)];\mathcal{R}\times SO(p))(f)=0$ for any
map germ $f:(\R^n,0)\lon (\R^p,0).$
We have the following general result as a simple corollary of Theorem 7.1 and Corollary 7.2.
\begin{Co}
Suppose that one of the conditions in Theorem 7.1 holds for $G.$
Then $\mathcal{M}(\mathcal{A}[G];\mathcal{R}\times G)(f)=0.$
\end{Co}
However, $\mathcal{M}(\mathcal{A}[SL(p,\R)];\mathcal{R}\times SL(p,\R)))(f)$
is an infinite dimensional vector space.
Following the examples in \S 6, there appear functional moduli for the classifications by $G$-congruence.
In the previous theory of singularities, functional moduli are usually unwelcome.
However, these give important information such as curvatures of curves etc.
\par
On the other hand, we have other relative moduli spaces for Lie subgroups $H<G<GL(p,\R)$:
\par
(1) $\displaystyle{\mathcal{M}(\mathcal{A}[G];\mathcal{A}[H])(f)=
\frac{T\mathcal{A}[G](f)}{T\mathcal{A}[H](f)}},$
\par
(2) $\displaystyle{\mathcal{M}(\mathcal{R}\times G;\mathcal{R}\times H)(f)=
\frac{T(\mathcal{R}\times G)(f)}{T(\mathcal{R}\times H)(f)}}.$
\par\noindent
Here, we call $\mathcal{M}(\mathcal{A}[G];\mathcal{A}[H])(f)$ a {\it relative infinitesimal moduli space of}
$f$ {\it with respect to $\mathcal{A}[G]$ and $\mathcal{A}[H].$}
We also call $\mathcal{M}(\mathcal{R}\times G;\mathcal{R}\times H)(f)$ a {\it relative infinitesimal moduli space of}
$f$ {\it with respect to $\mathcal{R}\times G$ and $\mathcal{R}\times H.$}
For $\mathcal{M}(\mathcal{A}[G];\mathcal{A}[H])(f),$ we have the
following exact sequence as $\R$-vector spaces:
\[
0\lon \frac{tf(\mathfrak{M}_n\theta (n))\cap \omega f_{[G]}(\theta [G]_0)}{tf(\mathfrak{M}_n\theta (n))\cap
\omega f_{[H]}(\theta [H]_0))}\lon \frac{\omega f_{[G]}(\theta [G]_0) }{\omega f_{[H]}(\theta [H]_0) }\stackrel{P}{\lon} \mathcal{M}(\mathcal{A}[G];\mathcal{A}[H])(f)\lon 0,
\]
where $P$ is defined by $P([\eta\circ f])=\{\eta\circ f\}.$
Since $T(\mathcal{A}[G])(f)=tf(\mathfrak{M}_n\theta (n))+\omega f_{[G]}(\theta [G]_0),$
$P$ is surjective and the kernel of $P$ is
\[
\frac{(tf(\mathfrak{M}_n\theta (n))+\omega f_{[H]}(\theta [H]_0))\cap \omega f_{[G]}(\theta [G]_0)}{\omega f_{[H]}(\theta [H]_0)}\cong
\frac{tf(\mathfrak{M}_n\theta (n))\cap \omega f_{[G]}(\theta [G]_0)}{tf(\mathfrak{M}_n\theta (n))\cap
\omega f_{[H]}(\theta [H]_0))}.
\]
It follows that
\[
\dim_\R \mathcal{M}(\mathcal{A}[G];\mathcal{A}[H])(f)\leq \dim_\R \frac{\omega f_{[G]}(\theta [G]_0) }{\omega f_{[H]}(\theta [H]_0) }\leq \dim _\R\frac{\theta [G]_0}{\theta [H]_0}.
\]
If $G=GL(p,\R)$ and $H=SL(p,\R)$, then
$\theta [GL(p,\R)]_0=\mathfrak{M}_p\theta (p)$ and 
$
\theta [SL(p,\R)]_0(p)=\ker ({\rm div}),
$
where ${\rm div} : \mathfrak{M}_p\theta (p)\lon \mathcal{E}_p$ is defined
by
\[{\rm div}\, \eta =\sum_{i=1}^p \frac{\partial \eta_i}{\partial y_i},\ \mbox{for}\ \eta=\sum_{i=1}^p \eta_i \frac{\partial }{\partial y_i}.
\]
The detailed calculations of
$\mathcal{M}(\mathcal{A}[GL(p,\R)];\mathcal{A}[SL(p,\R)])(f)$
were given in \cite{D-R}.
\par
For $\mathcal{M}(\mathcal{R}\times G;\mathcal{R}\times H)(f),$
we also have the following exact sequence:
\[
0\lon \frac{tf(\mathfrak{M}_n\theta (n))\cap \mathfrak{g}(f)}{tf(\mathfrak{M}_n\theta (n))\cap \mathfrak{h}(f)}\lon \frac{\mathfrak{g}(f)}{\mathfrak{h}(f)}\stackrel{\widetilde{P}}{\lon} \mathcal{M}(\mathcal{R}\times G;\mathcal{R}\times H)(f)\lon 0,
\]
where $\widetilde{P}$ is defined by $\widetilde{P}([X.f])=\{X.f\}.$
Since $T(\mathcal{R}\times G)(f)=tf(\mathfrak{M}_n\theta (n))+\mathfrak{g}(f),$
$\widetilde{P}$ is surjective and the kernel of $\widetilde{P}$ is
\[
\frac{(tf(\mathfrak{M}_n\theta (n))+\mathfrak{h}(f))\cap \mathfrak{g}(f)}{\mathfrak{h}(f)}
\cong
\frac{tf(\mathfrak{M}_n\theta (n))\cap \mathfrak{g}(f)}{tf(\mathfrak{M}_n\theta (n))\cap \mathfrak{h}(f)}.
\]
Moreover, we define the annihilator $\mathfrak{g}_f=\{X\in \mathfrak{g}\ |\ X.f=0\}.$
Then we have $\mathfrak{g}(f)\cong \mathfrak{g}/\mathfrak{g}_f.$
Since $\mathfrak{g}_f=M_p(\R)_f\cap \mathfrak{g},$ we have
\[
\frac{\mathfrak{g}(f)}{\mathfrak{h}(f)}\cong \frac{\mathfrak{g}/\mathfrak{g}_f}{\mathfrak{h}/\mathfrak{h}_f}
\cong \frac{\mathfrak{g}+M_p(\R)_f}{\mathfrak{h}+M_p(\R)_f}.
\]
Therefore, we have
\[
\dim_\R \mathcal{M}(\mathcal{R}\times G;\mathcal{R}\times H)(f)\leq
\dim_\R \frac{\mathfrak{g}+M_p(\R)_f}{\mathfrak{h}+M_p(\R)_f}\leq \dim_\R \frac{\mathfrak{g}}{\mathfrak{h}}=\dim G-\dim H <\infty.
\]
Since $\dim_\R GL(p,\R)-\dim_\R SL(p,\R)=1,$
$\dim_\R \mathcal{M}(\mathcal{R}\times GL(p,\R);\mathcal{R}\times SL(p,\R))(f)\leq 1.$
This means that the full affine geometry and the equi-affine geometry of map germs are
not so different.
If one of the conditions in Theorem 7.1 holds for $G$, then
\[
\mathcal{M}(\mathcal{A}[G];\mathcal{A}[H])(f)=\mathcal{M}(\mathcal{R}\times G;\mathcal{R}\times H)(f).
\]
It follows that 
\[
\dim_\R \mathcal{M}(\mathcal{A}[SO(4)];\mathcal{A}[\{1\}\oplus SO(3)])(f)\leq \dim SO(4)-\dim SO(3)
=3.
\]
\par
On the other hand, there is an interesting problem related to
$\mathcal{A}[G]$-equivalence and $\mathcal{R}\times G$-equivalence which we should investigate.
In \cite{M-S,M-S2,Saji} the normal forms with respect to
$SO(3)$-congruence of map germs $(\R^2,0)\lon (\R^3,0)$, which are $\mathcal{A}$-equivalent to the cuspidal edge or the swallowtail, are detected. Since there are no finitely determined map germs relative to $\mathcal{R}\times SO(3),$ these are not exact normal forms by the classification with respect
to $\mathcal{R}\times SO(3)=\mathcal{A}_0[SO(3)](2,3).$ They only give the Taylor polynomials of relatively lower orders by using ad hoc methods.
However, such normal forms give important geometric information of singular surfaces in $\R^3$
which are $\mathcal{A}$-equivalent to the cuspidal edge or the swallowtail.
All basic geometric invariants (i.e. various kinds of curvatures) at the origin are given by the
coefficients of these normal forms.
Therefore, we propose the following important but ambiguous problem:
\vskip3pt 
\par\noindent
{\bf Problem; semi-finite determinacy of map germs}: How can we determine the order of the Taylor polynomials
whose coefficients provide enough (or, all) geometric invariants with respect to
$\mathcal{R}\times G$ or $\mathcal{A}[G]$?
\vskip3pt
\par\noindent
For attacking this problem, we need extra new ideas beyond the Mather theory of $\mathcal{A}$-equivalence.
We suppose that the algebraic structure of the tangent space of the $\mathcal{A}[G]$-equivalence class
is one of the guideposts for solving the above problem.

\newpage
\begin{flushleft}
\textsc{
Shyuichi Izumiya
\\ Department of Mathematics
\\ Hokkaido University
\\ Sapporo 060-0810, Japan} 
\par
e-mail: izumiya@math.sci.hokudai.ac.jp

\bigskip
\textsc{
Masatomo Takahashi
\\ Muroran Institute of Technology,
\\ Muroran 050-8585, Japan}
\par
e-mail: masatomo@mmm.muroran-it.ac.jp

\bigskip
\textsc{
Hiroshi Teramoto
\\ Molecule \& Life Nonlinear Science Laboratory,
\\ Research Institute for Electronic Science,
\\ Hokkaido University
\\ Sapporo 001-0020, Japan} 
\par
e-mail: teramoto@es.hokudai.ac.jp

\end{flushleft}

\end{document}